\documentclass[12pt,reqno]{amsart}

\usepackage{mathabx}
\usepackage{amssymb}
\usepackage{bbm}
\usepackage{mathrsfs}
\usepackage[all]{xy}
\usepackage{float}

\newtheorem{facs}{Facts}[section]
\newtheorem{lem}[facs]{Lemma}
\newtheorem{prop}[facs]{Proposition}
\newtheorem{theo}[facs]{Theorem}

\newtheorem{facso}{Facts}[subsection]

\newtheorem{lemo}[facso]{Lemma}
\newtheorem{propo}[facso]{Proposition}
\newtheorem{coroo}[facso]{Corollary}
\newtheorem{theoo}[facso]{Theorem}

\theoremstyle{definition}

\newtheorem{ttt}[facs]{}

\newtheorem{defio}[facso]{Definition}

\newtheorem{algoo}[facso]{Algorithm}

\theoremstyle{remark}
\newtheorem{rem}[facs]{Remark}

\newtheorem{remo}[facso]{Remark}
\newtheorem{remso}[facso]{Remarks}
\newtheorem{exo}[facso]{Example}

\newcommand{\br}{ }
\newcommand{\brr}{, }

\newcommand{\Gal}{\mathop{\text{\rm Gal}}\nolimits}
\newcommand{\Pic}{\mathop{\text{\rm Pic}}\nolimits}
\newcommand{\NS}{\mathop{\text{\rm NS}}\nolimits}

\newcommand{\Frob}{\mathop{\text{\rm Frob}}\nolimits}
\newcommand{\bFrob}{\mathop{\text{\bf Frob}}\nolimits}

\renewcommand{\O}{\mathop{\text{\rm O}}\nolimits}
\newcommand{\SO}{\mathop{\text{\rm SO}}\nolimits}
\newcommand{\GO}{\mathop{\text{\rm GO}}\nolimits}

\newcommand{\End}{\mathop{\text{\rm End}}\nolimits}

\newcommand{\cha}{\mathop{\text{\rm char}}\nolimits}

\newcommand{\Tr}{\mathop{\text{\rm Tr}}\nolimits}
\newcommand{\rk}{\mathop{\text{\rm rk}}\nolimits}
\newcommand{\im}{\mathop{\text{\rm im}}\nolimits}
\newcommand{\id}{\mathop{\text{\rm id}}\nolimits}

\newcommand{\et}{\text{\rm \'et}}

\newcommand{\alg}{\text{\rm alg}}

\newcommand{\triv}{\text{\rm triv}}
\newcommand{\maxi}{{\!\text{\rm max}}}

\newcommand{\tors}{\text{\rm tors}}
\newcommand{\Hs}{\text{\rm Hs}}

\newcommand{\nr}{\text{\rm nr}}
\newcommand{\jump}{\text{\rm jump}}

\newcommand{\bbC}{{\mathbbm C}}
\newcommand{\bbF}{{\mathbbm F}}
\newcommand{\bbN}{{\mathbbm N}}
\newcommand{\bbQ}{{\mathbbm Q}}
\newcommand{\bbZ}{{\mathbbm Z}}

\newcommand{\calO}{{\mathscr{O}}}
\newcommand{\calS}{{\mathscr{S}}}

\newcommand{\frakp}{{\mathfrak p}}
\newcommand{\frakq}{{\mathfrak q}}

\newcommand{\frakz}{{\mathfrak z}}

\newcommand{\Pb}{{\text{\bf P}}}

\newcommand{\Exterior}{\mathchoice{{\textstyle\bigwedge}}%
 {{\bigwedge}}%
 {{\textstyle\wedge}}%
 {{\scriptstyle\wedge}}}

\newcounter{abc}
\newenvironment{abc}{\begin{list}{\rm \alph{abc}) }%
{\usecounter{abc} \leftmargin=0.0pt \labelsep=0.0pt %
\listparindent=0.0pt \labelwidth=0.0pt \parsep=\smallskipamount%
 \itemsep=0.0pt \topsep=0.0pt \partopsep=\smallskipamount}}{\end{list}}

\newcounter{iii}
\newenvironment{iii}{\begin{list}{\rm \roman{iii}) }%
{\usecounter{iii} \leftmargin=0.0pt \labelsep=0.0pt %
\listparindent=0.0pt \labelwidth=0.0pt \parsep=\smallskipamount%
 \itemsep=0.0pt \topsep=0.0pt \partopsep=\smallskipamount}}{\end{list}}

\let\hrighttoleftarrow=\righttoleftarrow
\renewcommand{\righttoleftarrow}{\!\hrighttoleftarrow}

\def\rightend#1#2{{%
 \leavevmode\nobreak\hskip .5em plus 1fil
 \penalty600 \hskip 0pt plus -1filll
 \vadjust{}\nobreak\hskip 0pt plus 1filll%
 #1\parfillskip=#2\relax \par}}

\def\eop{\ifmmode\rule[-22pt]{0pt}{1pt}\ifinner\tag*{$\square$}\else\eqno{\square}\fi\else\rightend{$\square$}{0pt}\fi}

\thanks{{\em Acknowledgements.} The authors would like to thank J.~Voight (Dartmouth College), who alerted them of an error in early version of this article, and the anonymous referees, whose comments turned out to be extremely helpful.}

\title[On the distribution of the Picard ranks]{On the distribution of the Picard ranks \\of the reductions of a
{\boldmath$K3$}~surface}

\begin{document}

\renewcommand{\thefootnote}{\fnsymbol{footnote}}
\author[Edgar Costa]{Edgar Costa${}^*$}

\footnotetext[1]{The first author was partially supported by the Simons Collaboration Grant \#550033.}

\address{Department of Mathematics \\ 77 Massachusetts Ave. \\ Bldg.\ 2-252B \\ Cambridge, MA 02139, USA}
\email{edgarc@mit.edu}
\urladdr{https://math.mit.edu/~edgarc/}

\author{Andreas-Stephan Elsenhans}

\address{Institut f\"ur Mathematik\\ Universit\"at W\"urzburg\\ Emil-Fischer-Stra\ss e 30\\ D-97074 W\"urzburg\\ Germany}
\email{stephan.elsenhans@mathematik.uni-wuerzburg.de}
\urladdr{https://www.mathematik.uni-wuerzburg.de/computeralgebra/team/elsenhans-step\discretionary{}{}{}han-\discretionary{}{}{}prof-dr/}

\author{J\"org Jahnel}

\address{\mbox{Department Mathematik\\ Univ.\ \!Siegen\\ \!Walter-Flex-Str.\ \!3\\ D-57068 \!Siegen\\ \!Germany}}
\email{jahnel@mathematik.uni-siegen.de}
\urladdr{http://www.uni-math.gwdg.de/jahnel}


\date{June~1,~2020}



\begin{abstract}
We report on our results concerning the distribution of the geometric Picard ranks of
$K3$~surfaces
under reduction modulo various primes. In the situation that
$\smash{\rk \Pic S_{\overline{K}}}$
is even, we introduce a quadratic character, called the jump character, such that
$\smash{\rk \Pic S_{\overline\bbF_{\!\frakp}} > \rk \Pic S_{\overline{K}}}$
for all good primes at which the character evaluates to
$(-1)$.
\end{abstract}

\maketitle
\thispagestyle{empty}

\section{Introduction}

Let
$S$
be a
$K3$~surface
over a number
field~$K$.
It~is a well-known fact that the geometric Picard rank
of~$S$
may not decrease under reduction modulo a good
prime~$\frakp$
of~$S$.
I.e.,~one always~has
\begin{equation}
\label{rkun}
\rk \Pic S_{\overline\bbF_{\!\frakp}} \geq \rk \Pic S_{\overline{K}} \, .
\end{equation}
%
It would certainly be interesting to understand the sequence
$\smash{(\rk \Pic S_{\overline\bbF_{\!\frakp}})_\frakp}$,
or at least the set of {\em jump~primes\/}
$$\Pi_\jump(S) := \{\,\frakp \text{ prime of } K \mid \frakp \text{ good for } S, \rk \Pic S_{\overline\bbF_{\!\frakp}} > \rk \Pic S_{\overline{K}}\,\} \, ,$$
for a given~surface. In an ideal case, one would be able to give a precise reason why the geometric Picard rank jumps at a given good~prime.
There are two well-known such~reasons.

\begin{iii}
\item
According to the Tate conjecture~\cite{KM}, the left hand side is always~even. Thus,~in the case that
$\rk \Pic S_{\overline{K}}$
is~odd, inequality (\ref{rkun}) is always strict and every good prime is a jump prime.
\item
Generalising this, if
$S$
has real multiplication (RM) by an endomorphism
field~$E$
and the integer
$(22 - \rk \Pic S_{\overline{K}})/[E\!:\!\bbQ]$
is odd, then again every good prime is a jump prime \cite[Theorem~1(2)]{Ch14}.
\end{iii}
It is known due to F.~Charles \cite[Theorem~1]{Ch14} that these are the only cases in which every good prime is a jump prime.

In~this article, we describe a third reason for a prime to jump, the {\em jump character}. It~was observed experimentally by the first author together with Yu.~Tschinkel \cite{CT} that, in the even rank case, one seems to have
$\smash{\liminf\limits_{B\to\infty}} \,\gamma(S,B) \geq \frac12$,~for
$$\gamma(S,B) := \frac{\#\{\frakp \in \Pi_\jump(S) \mid |\frakp| \leq B\}}{\#\{|\frakp| \leq B\}} \, .$$
Theorem~A shows that this observation is indeed true, except for some corner cases, which are simple to describe.
Moreover,~$\Pi_\jump(S)$
contains an entirely regular subset of exact density one~half.


In some sense, this result is complementary to those of F.\ Charles~\cite{Ch14}. In~fact, \cite[Theorem~1]{Ch14} shows, except for the cases i) and~ii), that there are always infinitely many {\em non-jump\/} primes.\medskip

\noindent
{\bf Theorem~A (Theorem \ref{K3_fl}).}
{\em
Let\/~$K$
be a number field and\/
$S$
a\/
$K3$~surface
over\/~$K$.
Moreover,~let\/
$\frakp \subset \calO_K$
be a prime of good reduction and residue
characteristic\/~$\neq \!2$.

\begin{abc}
\item
Then~there are two quantities, the\/ {\em discriminant}
$\Delta_{H^2}(S)$
of~$S$
and the\/ {\em discriminant of the Picard representation\/} or\/ {\em algebraic part of the discriminant\/}
$\Delta_{\Pic}(S)$
of~$S$
(cf.\ \ref{Delta_Hi}.ii) and \ref{Delta_NS}.a.ii) for precise definitions), such that the following equations hold,
\begin{align*}
\det(\Frob_\frakp \colon H^2_\et(S_{\overline{K}}, \bbQ_l(1)) \righttoleftarrow) &\!=\! \left( \frac{\Delta_{H^2}(S)}\frakp \right) \;\;\hspace{2cm} {\text and}\;\; \\ \det(\Frob_\frakp \colon T \righttoleftarrow) &\!=\! \left( \frac{\Delta_{H^2}(S)\Delta_{\Pic}(S)}\frakp \right) .
\end{align*}
Here,
$T := H_\alg^\perp \subset H^2_\et(S_{\overline{K}}, \bbQ_l(1))$
denotes the transcendental part of the~cohomology, and
$(\frac.\frakp)$
is the quadratic residue symbol
modulo~$\frakp$~\cite[Chapter~V, \S3]{Ne}.
\item
If\/
$\rk \Pic S_{\overline{K}}$
is even~then
\begin{equation}
\label{jumpeffect}
\left(\frac{\Delta_{H^2}(S)\Delta_{\Pic}(S)}\frakp\right) = -1 \quad\Longrightarrow\quad \rk \Pic S_{\overline\bbF_{\!\frakp}} \geq \rk \Pic S_{\overline{K}} + 2 \, .\end{equation}
In~other words, if\/
$K(\sqrt{\Delta_{H^2}(S)\Delta_{\Pic}(S)})/K$
is indeed a quadratic extension~then
$$\{\,\frakp \mid \frakp \text{ \rm inert in } K(\sqrt{\Delta_{H^2}(S)\Delta_{\Pic}(S)})/K\,\} \subseteq \Pi_\jump(S) \,.$$
\end{abc}
}\medskip

The~quadratic character
$\tau_S$,
given by
$$\frakp \mapsto \left(\frac{\Delta_{H^2}(S)\Delta_{\Pic}(S)}\frakp\right)$$
might be called the {\em transcendental character\/} of the
$K3$
surface~$S$.
Nevertheless, having implication (\ref{jumpeffect}) in mind, we prefer to call it the {\em jump character\/}
of~$S$,
at least in the even rank~case.
It~may happen that the quadratic extension, and hence the jump character, are~trivial. We~provide particular surfaces of this kind, defined
over~$\bbQ$,
in~Examples~\ref{Kummer}.a) and~\ref{jctriv}. These~are the corner cases mentioned~above.

One might think about the jump character also as~follows. The
\mbox{$\bbQ_l$-vector}
space
$T$
is equipped with the non-degenerate cup product pairing and acted upon by
$\Gal(\overline{K}/K)$.
Moreover,~the action is orthogonal with respect to the pairing, so that one has a continuous group homomorphism
$$\tau\colon \Gal(\overline{K}/K) \longrightarrow \O(T) \,.$$
There~is, however, no reason for
$\im \tau$
to be contained
in~$\SO(T)$,
in general, so the group homomorphism
$\det \tau\colon \Gal(\overline{K}/K) \to \{1, -1\}$
is usually non-trivial.
Moreover,~$\det \tau$
turns out to be independent
of~$l$,
cf.\ Proposition~\ref{quad_ext}.b).
The~article thus in essence describes the effects of
$\det \tau$
being non-trivial.

\subsubsection*{The criterion for non-triviality}\leavevmode

\noindent
Due~to its construction, the jump character is unramified at every prime of good reduction, cf.~Corollary \ref{jchar_nontriv}.a). On~the other hand, it may ramify at bad~primes. We~show this to be always the case when the singular reduction is of the mildest possible~type.\medskip

\noindent
{\bf Theorem~B (Corollary \ref{jchar_nontriv}.b)).}
{\em
Let\/~$K$
be a number field
and\/~$S$
a\/
$K3$~surface
over\/~$K$.
Moreover,~let\/
$\frakp \subset \calO_K$
be a prime of residue characteristic different
from\/~$2$.
Suppose~that\/
$S$
has a regular, projective model\/
$\underline{S}$
over\/~$\calO_{K,\frakp}$,
the geometric fibre
$\underline{S}_{\overline\frakp}$
of which has exactly one singular point, and assume this to be an ordinary double~point.\smallskip

\noindent
Then~the jump~character\/
$\smash{\tau_S = \big(\frac{\Delta_{H^2}(S)\Delta_{\Pic}(S)}\cdot\big)}$
ramifies
at\/~$\frakp$.}\medskip

In~addition, we present algorithms to compute the two characters for a given
surface~$S$
over~$\bbQ$,
a deterministic one for
$\Delta_{H^2}(S)$
and a statistical one for the jump~character. We~can deterministically compute the jump character in two particular situations, which together cover many, but not all~examples. These~situations occur when Corollary \ref{jchar_nontriv}.b) applies to at least a single bad prime, and when
$\Pic S_{\overline\bbQ}$
is known as a Galois module to such an extent that
$\smash{\Exterior^\maxi (\NS(S_{\overline{K}}) \otimes_\bbZ \bbQ)}$
can be determined. When~neither of these circumstances occurs, the statistical algorithm may still be~used.

\subsubsection*{An application. Rational curves on even rank\/
$K3$~surfaces}\leavevmode

\noindent
It~has been a long standing conjecture that every
$K3$~surface
over an algebraically closed field contains infinitely many rational curves. In full generality, it has been settled only recently by X.~Chen, F.~Gounelas, and C.~Liedtke \cite{CGL}. As~an application of our results, we show that the existence of infinitely many rational curves may be obtained rather easily in the situation that the jump character is non-trivial and the surface is otherwise~generic.\medskip

\noindent
{\bf Theorem~C (Theorem \ref{ratcrv}).}
{\em
Let\/~$K$
be a number field and\/
$S$
a\/
$K3$~surface
over\/~$K$.
Assume that\/
$\rk \Pic S_{\overline{K}}$
is even, that
$S_{\overline{K}}$
has neither real nor complex multiplication, and that\/
$\Delta_{H^2}(S)\Delta_{\Pic}(S)$
is a non-square
in\/~$K$.\smallskip

\noindent
Then\/
$S_{\overline{K}}$
contains infinitely many rational~curves.}\vspace{-0.5mm}

\subsubsection*{Notation. {\rm i)} Abstract Algebra, Algebraic Geometry, and Algebraic Number Theory}\leavevmode

\noindent
For~$K$
a field, we write
$\overline{K}$
for a fixed algebraic closure
of~$K$.
When
$K$
is a number field then by
$\calO_K$
we denote the ring of algebraic integers
in~$K$.
For~$\frakp \subset \calO_K$
a prime ideal, we let
$K_\frakp$
be the completion of
$K$
with respect
to~$\frakp$,
$\calO_{K,\frakp}$
the ring of integers in
$K_\frakp$,
and
$\bbF_{\!\frakp}$
the residue~field.

\looseness-1
Moreover,~we always let
$\smash{\Frob \in \Gal(\overline\bbF_{\!q}/\bbF_{\!q})}$
be the {\em geometric Frobenius\/} automorphism
$x \mapsto x^{1/q}$,
cf.~\cite[(1.15)]{DeWI}. Similarly, when
$K$
is a number field and
$\frakp \subset \calO_K$
a prime ideal,
$\smash{\Frob_\frakp \in \Gal(\overline{K}/K)}$
denotes an arbitrary lift of the geometric Frobenius. We~write
$\Frob$
and~$\Frob_\frakp$,
too, for the automorphisms of schemes and their cohomology groups, induced by
$\Frob$,
respectively~$\Frob_\frakp$,
via~functoriality. Let~us note that
$\smash{\Frob\colon \Pb^N_{\bbF_{q}} \to \Pb^N_{\bbF_{q}}}$
maps the geometric point
$(x_0:\cdots:x_N)$
to~$(x_0^q:\cdots:x_N^q)$.

Perhaps~deviating from a certain standard, we say that a proper variety
$S$
over~$K$
has {\em good reduction\/}
at~$\frakp$
when there exists a proper model
of~$S$
over~$\calO_K$
that has good reduction
at~$\frakp$
in the usual sense.\smallskip

\noindent
ii)
{\em Quadratic extensions.}
Let~$K$
be a field of characteristic
$\neq\! 2$
and
$L/K$
an at most quadratic field~extension. Then,~according to Kummer theory, there exists a unique class
$\Delta_L \in K^*/(K^*)^2$
such that
$\smash{L = K(\sqrt{\mathstrut u})}$
for any
$u \in \Delta_L$.
In~this situation, we shall also write
$\smash{K(\sqrt{\Delta_L})}$
for~$\smash{K(\sqrt{\mathstrut u})}$.

Assume~that
$K$
is a number field and
$\frakp \subset \calO_K$
a prime ideal of residue
characteristic~$\neq\! 2$,
at which
$L/K$
is~unramified. Then~the quadratic residue symbol
$\smash{(\frac{u}\frakp)}$~\cite[Chapter~V, Proposition~(3.5)]{Ne}
is independent of the choice of a
\mbox{$\frakp$-adic}
unit~$u \in \Delta_L$.
We~will therefore write
$\smash{(\frac{\Delta_L}\frakp)}$
instead
of~$\smash{(\frac{u}\frakp)}$.\smallskip

\noindent
iii)
{\em Characters.}
By a {\em character,} we always mean a continuous homomorphism from a topological group to a discrete abelian~group. A~{\em quadratic character\/} is a character
to~$\{1,-1\}$.

We~often describe a quadratic character
$\smash{\chi\colon \Gal(\overline{K}/K) \to \{1,-1\}}$,
when~$K$
is a number field, in the~form
$$\textstyle \frakp \mapsto (\frac\Delta\frakp) \,,$$
or simply
$(\frac\Delta\cdot)$,
for
$\Delta \in K^*/(K^*)^2$.
This~is supposed to mean that
$\chi(\Frob_\frakp) = (\frac\Delta\frakp)$
for every prime ideal
$\frakp \subset \calO_K$
of residue characteristic
$\neq \!2$
of the kind that
$\Delta$
is representable by a
$\frakp$-adic~unit.
Note that,
$\{1,-1\}$
being abelian,
$\chi(\Frob_\frakp)$
is well defined for every quadratic character. Moreover,~the values at the Frobenii
$\Frob_\frakp$
determine
$\chi$
uniquely, due to the Chebotarev density~theorem.

When
$V$
is a one-dimensional vector space over a
field~$F$,
equipped with the discrete topology and acted upon continuously by a topological
group~$G$,
we denote by
$[V]\colon G \to F^*$
the character given~by
$[V](g) = a$,
for
$a \in F^*$
the scalar satisfying
$g \!\cdot\! v = av$
for every
$v \in V$.
This~notation is due to T.~Saito~\cite{Sa}.

\subsubsection*{Computations}\leavevmode

\noindent
All~computations are done using {\tt magma}~\cite{BCP}, {\tt sage}~\cite{St}, and {\tt C++}, including the libraries {\tt FLINT}~\cite{HJP} and {\tt NTL}~\cite{Sh}. For~point counting on the examples being quartic surfaces, we used the software developed by the first author, which is publicly available at {\tt https://github.com/edgarcosta/controlledreduction}.

\section{The jump character}

\subsection{The determinant of\/
$\bFrob$
and the relationship with the sign in the functional equation}\leavevmode

\noindent
Let~$S$
be a smooth, proper variety over a finite
field~$\bbF_{\!q}$
of
characteristic~$p > 0$.
Then~$\Frob$
acts linearly on the
\mbox{$l$-adic}
cohomology modules
$\smash{H^i_\et(S_{\overline\bbF_{\!q}}, \bbZ_l(j))}$.
The~characteristic polynomial
$\smash{\Phi_j^{(i)}}$
of~$\Frob$
is independent of the choice
of~$l \neq p$
and has in fact rational coefficients~\cite[Th\'eor\`eme~(1.6)]{DeWI}. In~particular, the determinant
of~$\Frob$
is a rational number and independent
of~$l \neq p$.

In this section, we discuss the behaviour of
$\det \Frob$.
The~theorem below, which is essentially a summary of results of P.\ Deligne and J.\ Suh, shows that each statement
on~$\det\Frob$
may be translated into a statement about the sign in the functional~equation.

\begin{theoo}[Deligne, Suh]
\begin{abc}
\item
The polynomial\/
$\smash{\Phi_j^{(i)} \in \bbQ[T]}$
fulfils the functional equa\-tion
\begin{equation}
\label{eq_eins}
T^N \Phi(q^{i-2j}/T) = \pm q^{\frac{N}2(i - 2j)} \Phi(T) \, ,
\end{equation}
for\/~$\smash{N := \rk H^i_\et(S_{\overline\bbF_{\!q}}, \bbZ_l(j))}$.
\item
The sign in the functional equation (\ref{eq_eins}) is that~of
$$\det(-\Frob\colon H^i_\et(S_{\overline\bbF_{\!q}}, \bbQ_l(j)) \righttoleftarrow) = (-1)^N \det(\Frob\colon H^i_\et(S_{\overline\bbF_{\!q}}, \bbQ_l(j)) \righttoleftarrow) \, .$$
This is a rational number, the sign of which is independent of the Tate twist, i.e.\ of the choice
of\/~$j$.
\item[{\rm c.i) }]
If\/~$i$
is even, then\/
$\smash{\det(-\Frob\colon H^i_\et(S_{\overline\bbF_{\!q}}, \bbQ_l(i/2)) \righttoleftarrow)}$
is either\/
$(+1)$
or\/~$(-1)$.
In~other words, the determinant gives the sign in (\ref{eq_eins})~exactly.
\item[{\rm ii) }]
If\/~$i$
is odd then\/
$N$~is
even and in (\ref{eq_eins}), the plus sign always holds.
\end{abc}
\medskip

\noindent
{\bf Proof.}
{\em
a) and~b)
Let~us write
$\Phi$
for~$\smash{\Phi_j^{(i)}}$.
The~polynomials on both sides of (\ref{eq_eins}) then have the same roots as,
with~$z$,
the~number
$\smash{\overline{z} = \frac{q^{i-2j}}z}$
is a root
of~$\Phi$,
too \cite[Corollaire~(3.3.9)]{DeWII}, and has the same multiplicity. To~show that they perfectly agree, let us adopt the convention that
$\Phi$
is~monic. Then~the leading coefficient of the polynomial on the left hand side is equal to the constant term
of~$\Phi$.
By~definition, this is the determinant of
$(-\Frob)$
on
$\smash{H^i_\et(S_{\overline\bbF_{\!q}}, \bbQ_l(j))}$,
which is known to be a rational number and of absolute value
$\smash{q^{\frac{N}2(i - 2j)}}$
\cite[Corollaire~(3.3.9)]{DeWII}.
Thus,~a) follows, together with the first assertion of~b). The~final claim is clear,~since
$$\det(-\Frob\colon H^i_\et(S_{\overline\bbF_{\!q}}, \bbQ_l(j)) \righttoleftarrow) = q^{-Nj} \!\cdot \det(-\Frob\colon H^i_\et(S_{\overline\bbF_{\!q}}, \bbQ_l) \righttoleftarrow) \, .$$
c.i)
As~$\det\Frob = (-1)^N \Phi(0)$
for
$\Phi$
monic, this can be read off the functional equa\-tion
$T^N \Phi(1/T) = \pm\Phi(T)$.\smallskip

\noindent
ii)
If~$S$
is projective then, by Poincar\'e duality and the hard Lefschetz theorem \cite[Th\'e\-o\-r\`eme (4.1.1)]{DeWII}, there is a non-degenerate~pairing
$$\smash{H^i_\et(S_{\overline\bbF_{\!q}}, \bbQ_l(j)) \times H^i_\et(S_{\overline\bbF_{\!q}}, \bbQ_l(j)) \to \bbQ_l(2j - i)}$$
that is compatible with the action
of~$\Frob$.
It~is, moreover, alternating since
$i$
is supposed to be odd, cf.\ \cite[Chapter~5, Section~6, \S11]{Sp}. The~assertion follows directly from this \cite[(2.6)]{DeWI}. Cf.~the remarks after~\cite[Corollaire~(4.1.5)]{DeWII}. The proper non-projective case has only recently been settled by J.~Suh \cite[Corollary~2.2.3 and Corollary~3.3.5]{Su}.%
}%
\eop
\end{theoo}

\begin{remo}
The~polynomials
$\smash{\Phi_j^{(i)} \in \bbQ[T]}$
occurring as characteristic polynomials
of~$\Frob$
have remarkable properties, which were established mainly by P.~Deligne, B.~Mazur, and A.~Ogus. Details are summarised in the article \cite{EJ15} of the second and third~authors. In~the proof above, the only property that was used is that every complex root
of\/~$\smash{\Phi_j^{(i)}}$
is of absolute
value\/~$q^{i/2 - j}$.
This~was first proven by P.~Deligne in~\cite[Th\'e\-o\-r\`eme~(1.6)]{DeWI} for the projective case and later in \cite[Corollaire~(3.3.9)]{DeWII}, in~general. The~assertion had been formulated by A.\,Weil as a part of his famous~conjectures.
\end{remo}

\subsection{The discriminant of the
{\boldmath $H^i$}-representation}\leavevmode

\noindent
Let~us start by recalling some facts on
\mbox{$l$-adic}
cohomology.

\begin{propo}
\label{quad_ext}
Let\/~$K$
be a field
and\/~$S$
a smooth and proper\/
\mbox{$K$-scheme}.

\begin{abc}
\item
Then, for all prime numbers\/
$l \neq \cha K$
and all integers\/
$i$
and\/~$j$,
associated with the one-dimensional\/
\mbox{$\bbQ_l$-vector}
space\/
$\Exterior^\maxi H^i_\et(S_{\overline{K}}, \bbQ_l(j))$,
there is the~character
$$[\det H^i(S_{\overline{K}}, \bbQ_l(j))] \colon \Gal(\overline{K}/K) \longrightarrow \bbQ_l^*$$
of the absolute Galois group
of\/~$K$.
\item
Suppose~that\/
$i$
is even and that\/
$S$
is pure of
dimension\/~$i$.
Then~the character\/
$\smash{[\det H^i(S_{\overline{K}}, \bbQ_l(i/2))]}$
has values
in\/~$\{1, -1\} \subset \bbQ_l^*$
and is independent
of\/~$l$.
\end{abc}\smallskip

\noindent
{\bf Proof.}
{\em
a)
This~follows from the functoriality of
\mbox{$l$-adic}
cohomology, together with the fact that every
$\smash{\sigma \in \Gal(\overline{K}/K)}$
induces an automorphism of schemes
of~$S_{\overline{K}}$.\smallskip

\noindent
b)
By~Poincar\'e duality~\cite[Exp.~XVIII, Th\'e\-o\-r\`eme~3.2.5]{SGA4}, there is a canonical non-degenerate~pairing
$\smash{s\colon H^i_\et(S_{\overline{K}}, \bbQ_l(i/2)) \times H^i_\et(S_{\overline{K}}, \bbQ_l(i/2)) \to \bbQ_l}$
that is compatible with the action
of~$\smash{\Gal(\overline{K}/K)}$.
Here,
$i$
is assumed even, so the pairing
$s$
is~symmetric. According~to a standard fact from linear algebra \cite[Def.~2.9]{War},
$s$~induces
another symmetric~pairing
$$\Exterior^\maxi(s)\colon \Exterior^\maxi H^i_\et(S_{\overline{K}}, \bbQ_l(i/2)) \times \Exterior^\maxi H^i_\et(S_{\overline{K}}, \bbQ_l(i/2)) \longrightarrow \bbQ_l$$
that is again non-degenerate. The~action
of~$\smash{\Gal(\overline{K}/K)}$
is orthogonal with respect
to~$s$,
which implies that the character
$\smash{[\det H^i(S_{\overline{K}}, \bbQ_l(i/2))]}$
must have values
in~$\{1, -1\} \subset \bbQ_l^*$.

For~the case that
$K$
is a number field, independence
of~$l$
is easily reduced to the Weil conjectures, proven by P.~Deligne \cite[Corollaire~(3.3.9)]{DeWII}, using the Chebotarev density theorem together with the smooth specialisation theorem for cohomology groups \cite[Exp.\ XVI, Corollaire 2.3]{SGA4}. In~general, the result has been established by T.~Saito~\cite[Corollary~3.3]{Sa}.

N.B. In the notation for the characters, we write
$\det$
instead
of~$\Exterior^\maxi$.
This~convention follows~\cite{Sa}.
}
\eop
\end{propo}

\begin{defio}
\label{Delta_Hi}
\begin{iii}
\item
In the situation of part~b), we denote by
$L_S$
the extension field of
$K$
that corresponds to
$\ker\, [\det H^i(S_{\overline{K}}, \bbQ_l(i/2))]$
under the Galois correspondence. By~construction,
$L_S/K$
is an at most quadratic~extension.
\item\looseness-1
If
$\cha K \neq 2$
then we denote the class
in~$K^*/(K^*)^2$
that yields the field extension
$L_S/K$
by~$\Delta_{H^i}(S)$
and call it the {\em discriminant of the
\mbox{$H^i$-representation\/}}
of~$S$.
\end{iii}
\end{defio}

\begin{lemo}
\label{unverzw_Hi}
Let\/~$K$
be a number field
and\/~$S$
a smooth and proper\/
\mbox{$K$-scheme}.
Moreover,~let\/
$\frakp \subset \calO_K$
be a prime at which\/
$S$~has
good~reduction. 

\begin{abc}
\item
If\/~$l$
is a prime different from the residue characteristic
of\/~$\frakp$
then, for every\/
$j \in \bbZ$,
the\/
\mbox{$\smash{\Gal(\overline{K}/K)}$-representation\/}
$\smash{H^i_\et(S_{\overline{K}}, \bbQ_l(j))}$
is unramified
at\/~$\frakp$.
\item
Suppose that\/
$S$
is pure of
dimension\/~$i$,
for an even
integer\/~$i$.
Then~the quadratic character\/
$\smash{\big(\frac{\Delta_{H^i}(S)}\cdot\big) = [\det H^i(S_{\overline{K}}, \bbQ_l(i/2))]}$
is unramified
at\/~$\frakp$.
Equivalently, the splitting field\/
$L_S$
is unramified
at\/~$\frakp$.
\end{abc}\smallskip

\noindent
{\bf Proof.}
{\em
a)
For this, it suffices to consider the restriction of the representation to the decomposition group
$\smash{D_\frakp \cong \Gal(\overline{K}_\frakp/K_\frakp)}$.
This~coincides with the natural action of
$\smash{\Gal(\overline{K}_\frakp/K_\frakp)}$
on
$\smash{H^i_\et(S_{\overline{K}_\frakp}, \bbQ_l(j)) \cong H^i_\et(S_{\overline{K}}, \bbQ_l(j))}$,
according to invariance of \'etale cohomology under extensions of separably closed fields \cite[Exp.~XII, Corollaire~5.4]{SGA4}.
Moreover,~by the smooth specialisation theorem for cohomology groups \cite[Exp.\ XVI, Corollaire 2.2]{SGA4},
$\smash{H^i_\et(S_{\overline{K}_\frakp}, \bbQ_l(j)) \cong H^i_\et(S_{\overline\bbF_{\!\frakp}}, \bbQ_l(j))}$,
which shows that this cohomology vector space is acted upon via the quotient
$\smash{\Gal(\overline\bbF_{\!\frakp}/\bbF_{\!\frakp}) \cong D_\frakp/I_\frakp}$.
I.e.,~the inertia group
$I_\frakp$
fixes
$\smash{H^i_\et(S_{\overline{K}}, \bbQ_l(j))}$
pointwise, as~required.\smallskip

\noindent
b)
This is a direct consequence of~a).
}
\eop
\end{lemo}

\begin{remso}
\begin{iii}
\item
In~particular, for every exponent
$k \in \bbN$,
the splitting field of
$\smash{H^i_\et(S_{\overline{K}}, \bbZ/l^k\bbZ(j))}$
is unramified at every prime ideal
$\frakp$
of good reduction.
\item
If~the residue characteristic
is\/~$\neq \!2$
then one~has
\begin{equation}
\label{Hi_char}
\det(\Frob_\frakp\colon H^i_\et(S_{\overline{K}}, \bbQ_l(i/2)) \righttoleftarrow) = \left( \frac{\Delta_{H^i}(S)}\frakp \right) \, .
\end{equation}
\item
For~any number field
$K' \supseteq K$,
one has
$\Delta_{H^i}(S_{K'}) = \Delta_{H^i}(S) \cdot (K'{}^*)^2 \in K'/(K'{}^*)^2$.
\item
The name ``discriminant'' has not been chosen at~random. Indeed, Proposition~\ref{quad_ext} allows a generalisation to families
$\pi\colon F \to X$
over a general base scheme. The~quadratic field extension then goes over into a twofold covering
$\varrho_\pi \colon Y \to X$,
ramified at most over the discriminant~locus. 
If~$\pi$
is sufficiently reasonable then
$\varrho_\pi$
is given by
$w^2 = \Delta$,
for
$\Delta$
a normalised version of the discriminant of the~family. Furthermore,~for every non-singular member
$S = F_x$,
the value
$\Delta(x)$
belongs to the
class~$\Delta_{H^i}(S)$.

We~plan to report about the interpretation of the quantities
$\Delta_{H^i}(S)$
as actual discriminants, as well as some applications thereof, in a forthcoming paper.
\end{iii}
\end{remso}

\subsection{Surfaces---The discriminant of the N\'eron--Severi representation}\leavevmode

\begin{propo}
Let\/~$K$
be a field
and\/~$S$
a smooth projective surface
over~\/$K$.

\begin{abc}
\item
Then, associated with the one-dimensional\/
\mbox{$\bbQ$-vector}
space\/
$\Exterior^\maxi (\NS(S_{\overline{K}}) \otimes_\bbZ \bbQ)$,
there is the~character
$$[\det (\NS(S_{\overline{K}}) \otimes_\bbZ \bbQ)] \colon \Gal(\overline{K}/K) \longrightarrow \bbQ^*$$
of the absolute Galois group
of\/~$K$.
\item
The character\/
$[\det (\NS(S_{\overline{K}}) \otimes_\bbZ \bbQ)]$
takes values only in\/
$\{1, -1\} \subset \bbQ^*$.
\end{abc}\smallskip

\noindent
{\bf Proof.}
{\em
a)
This~follows from the functoriality of the N\'eron--Severi group, together with the fact that every
$\smash{\sigma \in \Gal(\overline{K}/K)}$
induces an automorphism of schemes
of~$S_{\overline{K}}$.\smallskip

\noindent
b)
Every~$\sigma \in \Gal(\overline{K}/K)$
induces an automorphism of the N\'eron--Severi group
$\NS(S_{\overline{K}})$,
in particular one of
$\NS(S_{\overline{K}})_\tors$,
and consequently one of the torsion-free
\mbox{$\bbZ$-module}
$\NS(S_{\overline{K}})/\NS(S_{\overline{K}})_\tors$.
As~that is a full rank lattice
in~$\NS(S_{\overline{K}}) \otimes_\bbZ \bbQ$,
it induces a lattice in the one-dimensional vector space
$\smash{\Exterior^\maxi (\NS(S_{\overline{K}}) \otimes_\bbZ \bbQ)}$,
which~must be respected by the action
of~$\smash{\Gal(\overline{K}/K)}$.
The~assertion immediately follows from this.
}
\eop
\end{propo}

\begin{defio}
\label{Delta_NS}
\begin{abc}
\item[]
\begin{iii}
\item[a.i) ]
In the situation of part~b), we denote by
$L_{S,\alg}$
the extension field of
$K$
that corresponds to
$\ker\, [\det (\NS(S_{\overline{K}}) \otimes_\bbZ \bbQ)]$
under the Galois correspondence.
\item[ii) ]\looseness-1
If~$\cha K \neq 2$
then we denote the class
in~$K^*/(K^*)^2$
that yields the field extension
$\smash{L_{S,\alg}/K}$
by~$\Delta_{\NS}(S)$
and call it the {\em discriminant of the N\'eron--Severi rep\-re\-sen\-ta\-tion} or the {\em algebraic part of the discriminant}
of~$S$.
\end{iii}
\item[b) ]
For~surfaces such that\/
$H^1(S_{\overline{K}}, \bbQ_l) = 0$,
one has
$\Pic(S_{\overline{K}}) \otimes_\bbZ \bbQ \cong \NS(S_{\overline{K}}) \otimes_\bbZ \bbQ$.
In~this case, one may write
$\Delta_{\Pic}(S)$
instead of
$\Delta_{\NS}(S)$
and speak of the {\em discriminant of the Picard rep\-re\-sen\-ta\-tion}. Similarly,~let us then write
$[\det (\Pic(S_{\overline{K}}) \otimes_\bbZ \bbQ)]$
instead of
$[\det (\NS(S_{\overline{K}}) \otimes_\bbZ \bbQ)]$.
\end{abc}
\end{defio}

\begin{remo}
The~algebraic part of the discriminant
of~$S$
should not be confused with the discriminants of
$\NS(S)$
and~$\NS(S_{\overline{K}})$
as~lattices. Instead,~one might think about it as~follows.

Since
$\NS(S_{\overline{K}}) \otimes_\bbZ \bbQ$
is a finite-dimensional
\mbox{$\bbQ$-vector}
space, there is a smallest finite field extension
of~$K$,
over which all elements
of~$\NS(S_{\overline{K}}) \otimes_\bbZ \bbQ$
are defined, the splitting
field~$L$
of~$\NS(S_{\overline{K}}) \otimes_\bbZ \bbQ$.
Then~$\smash{\Gal(\overline{K}/K)}$
acts on
$\NS(S_{\overline{K}}) \otimes_\bbZ \bbQ$
via its quotient
$\Gal(L/K)$
and the action of this quotient is~faithful.

Therefore,~the action of
$\smash{\Gal(\overline{K}/K)}$
on
$\smash{\Exterior^\maxi (\NS(S_{\overline{K}}}) \otimes_\bbZ \bbQ)$
factors via
$\Gal(L/K)$,
too.
On~the other hand,
$\smash{\Exterior^\maxi (\NS(S_{\overline{K}}}) \otimes_\bbZ \bbQ)$
is one-dimensional, so the action of the finite group
$\Gal(L/K)$
must factor
via~$(\bbQ^*)_\tors = \{1,-1\}$.
The~stabiliser is of index at
most~$2$
in~$\smash{\Gal(\overline{K}/K)}$
and hence the splitting field of
$\smash{\Exterior^\maxi (\NS(S_{\overline{K}}}) \otimes_\bbZ \bbQ)$
is an at most quadratic extension
$\smash{K(\sqrt{\Delta_{\NS}(S)}) = L_{S,\alg} \subseteq L}$.
\end{remo}

\begin{lemo}\looseness-1
\label{unverzw}
Let\/~$K$
be a number field
and\/~$S$
a smooth projective surface
over\/~$K$.
Moreover,~let\/
$\frakp \subset \calO_K$
be a prime at which\/
$S$~has
good~reduction.

\begin{abc}
\item
Then~the splitting field of\/
$\NS(S_{\overline{K}}) \otimes_\bbZ \bbQ$
is unramified
at\/~$\frakp$.
\item
The~character\/
$\smash{\big(\frac{\Delta_{\NS}(S)}\cdot\big) = [\det \NS(S_{\overline{K}}) \otimes_\bbZ \bbQ]}$
is unramified
at\/~$\frakp$.
Equivalently, the splitting field\/
$L_{S,\alg}$
is unramified
at\/~$\frakp$.
\end{abc}\smallskip

\noindent
{\bf Proof.}
{\em
a)
The first Chern class homomorphism factors via the N\'eron--Severi group, i.e.\ via algebraic equivalence,
$$
\xymatrix{
\Pic(S_{\overline{K}}) \otimes_\bbZ \bbQ_l \ar@{->}[rr]^{c_1} \ar@{->>}[d] && H^2_\et(S_{\overline{K}}, \bbQ_l(1)) \\
\NS(S_{\overline{K}}) \otimes_\bbZ \bbQ_l \ar@{->}[rru] \,. &&
}
$$
Indeed,~$c_1$
factors via numerical equivalence, since the intersection pairing on 
$\Pic(S_{\overline{K}}) \otimes_\bbZ \bbQ_l$
is compatible with the cup product pairing
on~$H^2_\et(S_{\overline{K}}, \bbQ_l(1))$.
Moreover, Matsusaka's theorem \cite[Theorem~4]{Ma}, cf.~\cite[paragraph~3.2.7]{An}, shows that algebraic equivalence coincides with numerical equivalence, already on
$\Pic(S_{\overline{K}}) \otimes_\bbZ \bbQ$.

Now~write
$L$
for the splitting field of
$\smash{\NS(S_{\overline{K}}) \otimes_\bbZ \bbQ}$
and assume~that
$L$
would ramify
at~$\frakp$.
By~definition,
$\Gal(L/K)$
acts faith\-fully
on~$\smash{\NS(S_{\overline{K}}) \otimes_\bbZ \bbQ}$.
Choose~a prime
$\frakq$
of~$L$
lying
above~$\frakp$
and a non-trivial element
$\sigma \in \Gal(L_\frakq/K_\frakp^{\rm n}) \subseteq \Gal(L_\frakq/K_\frakp)$,
for
$K_\frakp^{\rm n}$
the maximal unramified subfield
of~$L_\frakq$.
Then~$\sigma$
acts non-trivially on the image of the first Chern class homomorphism
$$\smash{c_1\colon \NS(S_{\overline{K}}) \otimes_\bbZ \bbQ_l \hookrightarrow H^2_\et(S_{\overline{K}}, \bbQ_l(1)) \, .}$$
This,~however, is in contradiction with the smooth specialisation theorem for co\-homol\-ogy groups, as seen~before.

\smallskip

\noindent
b)
This follows immediately from~a).
}
\eop
\end{lemo}

\begin{remso}
\begin{iii}
\item
If~the residue characteristic is
not\/~$2$
then one~has
\begin{equation}
\label{det_NS}
\det(\Frob_\frakp\colon (\NS(S_{\overline{K}}) \otimes_\bbZ \bbQ) \righttoleftarrow) = \left( \frac{\Delta_{\NS}(S)}\frakp \right) \, .
\end{equation}
\item
If
$K' \supseteq K$
is a number field
extending~$K$,
then one~has
$$\Delta_{\NS}(S_{K'}) = \Delta_{\NS}(S) \cdot (K'{}^*)^2 \in K'/(K'{}^*)^2.$$
\item
It is worthwhile to observe that, for
$K3$~surfaces,
the assertion of Lemma~\ref{unverzw} is still true when
$S$~has
bad reduction
at~$\frakp$
of the mildest possible~form. Cf.~Corollary \ref{Pic_split},~below.\smallskip
\end{iii}
\end{remso}

\subsection{{\boldmath$K3$}~surfaces}\leavevmode

\noindent
There~is a strong relation, which is established for
$K3$~surfaces,
but relies on Tate's and Serre's conjectures in general, between the Galois action on
\mbox{$l$-adic}
cohomology and the variation of the geometric Picard ranks under reduction modulo various primes~\cite{EJ12b,CT}. From our point of view, this is in fact the main application of the constructions presented so~far.

\begin{facso}
\label{K3_gen}
Let\/~$S$
be a\/
$K3$~surface
over a base
field\/~$K$.

\begin{abc}
\item
Then\/~$\Pic S_{\overline{K}}$
is a free abelian group of rank at
most\/~$22$.
If\/~$K$
is of characteristic zero then the rank is at
most\/~$20$.
If\/~$K$
is finite then the rank is~even.
\item
If\/~$K$
is finite then\/
$\rk \Pic S_{\overline{K}}$
is equal to the number [counted with multiplicities] of all eigenvalues
of\/~$\Frob$
on\/
$H^2_\et(S_{\overline{K}}, \bbQ_l(1))$
that are roots of~unity.
\end{abc}\medskip

\noindent
{\bf Proof.}
{\em
a)
The~first statement is found, e.g., in \cite[Chapter~17, first formula of Section~2]{Hu}. The second one is \cite[Chapter~17, formula~(1.1)]{Hu} in the case that
$K=\bbC$,
while the general case follows from this in view of \cite[Chapter~17, Lemma~2.2]{Hu}. The final claim is a direct consequence of~b).\smallskip

\noindent
b)
See \cite[Chapter~17, Corollary~2.9 and the arguments given before]{Hu}. Note that this result is an application of the Tate conjecture, which has been shown for
$K3$~surfaces
over finite fields by to the combined work of several people, most notably F.~Charles \cite{Ch13}, M.~Lieb\-lich, D.\ Maulik, and A.\ Snowden~\cite{LMS}, K.~Madapusi Pera \cite{MP}, as well as W.~Kim and K.\ Madapusi Pera \cite{KM}.
}
\eop
\end{facso}

For~$K$
an arbitrary field and
$l \neq \cha K$
a prime number, there is a canonical orthogonal decomposition
\begin{equation}
\label{split}
H^2_\et(S_{\overline{K}}, \bbQ_l(1)) = H_\alg \oplus T \, .
\end{equation}
Here,~$H_\alg = c_1(\Pic(S_{\overline{K}}) \otimes_\bbZ \bbQ_l)$
is clearly
\mbox{$\Gal(\overline{K}/K)$-invariant}.
Moreover,~$T := H_\alg^\perp$
is
\mbox{$\smash{\Gal(\overline{K}/K)}$-invariant},
too, as the Galois action is~orthogonal.

In~the particular case that
$K$
is a number field, let
$\frakp \subset \calO_K$
be any prime of good reduction and of residue characteristic different
from~$l$.
Then~$\smash{\Frob_\frakp \in \Gal(\overline{K}/K)}$
is determined only up to conjugation. But this suffices to have well-defined eigenvalues and a well-defined determinant
of~$\Frob_\frakp$,
associated with any vector space being acted upon
by~$\smash{\Gal(\overline{K}/K)}$.
In~particular,
\begin{equation}
\label{det_Formel}
\det(\Frob_\frakp \colon T \righttoleftarrow) = \frac{\det(\Frob_\frakp \colon H^2_\et(S_{\overline{K}}, \bbQ_l(1)) \righttoleftarrow)}{\det(\Frob_\frakp \colon \Pic S_{\overline{K}} \righttoleftarrow)} \, .
\end{equation}
Our main theoretical observation on the distribution of the Picard ranks of the reductions is then as~follows.

\begin{propo}[Rank jumps]
\label{jumps}
Let\/~$S$
be a\/
$K3$~surface
over a number
field\/~$K$
and\/
$\frakp \subset \calO_K$
a prime of good reduction.
Assume that\/
$\rk \Pic S_{\overline{K}}$
is~even. Then the following is~true:\smallskip

\noindent
If\/
$\det(\Frob_\frakp \colon T \righttoleftarrow) = -1$,
then\/
$\rk \Pic S_{\overline\bbF_{\!\frakp}} \geq \rk \Pic S_{\overline{K}} + 2$.\smallskip

\noindent
{\bf Proof.}
{\em
Choose a prime number
$l$,
different from the residue \mbox{characteristic
of~$\frakp$.}
Then~one has
$\smash{H^2_\et(S_{\overline{K}}, \bbQ_l(1)) \cong H^2_\et(S_{\overline{K}_\frakp}, \bbQ_l(1))}$
\cite[Exp.~XII, Corollaire~5.4]{SGA4}, as well as
$\smash{H^2_\et(S_{\overline{K}_\frakp}, \bbQ_l(1)) \cong H^2_\et(S_{\overline\bbF_{\!\frakp}}, \bbQ_l(1))}$
\cite[Exp.~XVI, Corollaire~2.2]{SGA4}.
Consequently, by transport of structure, the orthogonal decomposition (\ref{split}) carries over into
\begin{equation}
\label{split_p}
H^2_\et(S_{\overline\bbF_{\!\frakp}}, \bbQ_l(1)) = H_\alg \oplus T \, .
\end{equation}
Note~that, as a consequence of its construction,
$T$
may well contain algebraic~classes.

Moreover,~under the first isomorphism, the action
of~$\smash{\Gal(\overline{K}_\frakp/K_\frakp)}$
is compatible with that of the decomposition group
$\smash{D_\frakp \subset \Gal(\overline{K}/K)}$,
while the second isomorphism shows that
$\smash{\Gal(\overline{K}_\frakp/K_\frakp)}$
acts via its quotient
$\smash{\Gal(K^\nr_\frakp/K_\frakp)}$.
In~particular, the action of
$\smash{\Frob \in \Gal(\overline\bbF_{\!\frakp}/\bbF_{\!\frakp})}$
on
$\smash{H^2_\et(S_{\overline\bbF_{\!\frakp}}, \bbQ_l(1))}$
agrees with that of any
$\smash{\Frob_\frakp \in \Gal(\overline{K}/K)}$
on~$\smash{H^2_\et(S_{\overline{K}}, \bbQ_l(1))}$.
For~instance,
$\Frob$
and
$\Frob_\frakp$
have the same eigenvalues on
$H_\alg$,
as well as
on~$T$.

By~Lemma~\ref{unverzw}.a), the splitting field
of~$\Pic S_{\overline{K}}$
is a number field unramified
at~$\frakp$.
Therefore,~$\smash{\Gal(K^\nr_\frakp/K_\frakp)}$
acts
on~$H_\alg$
via a finite quotient~group. In~particular, there exists an integer
$e>0$
such that
$\Frob^e$
acts~trivially. Consequently,~all eigenvalues of
$\Frob$
on~$H_\alg$
are roots of unity.

In~view of Fact~\ref{K3_gen}.b), we need to show that
$\Frob$
acts
on~$T$
with at least two eigenvalues being roots of~unity. For~this, let us observe that each eigenvalue is of absolute
value~$1$,
so that those different from
$1$
and~$(-1)$
come in pairs
$\{z, \overline{z}\}$
of complex conjugates.
As~$z\overline{z} = 1$
and
$\det(\Frob_\frakp \colon T \righttoleftarrow) = -1$,
one of the eigenvalues must be equal
to~$(-1)$.
Finally,~as
$\dim T = 22 - \rk \Pic S_{\overline{K}}$
is even, a further
eigenvalue~$1$
is~enforced. This~completes the~proof.
}
\eop
\end{propo}

\begin{remso}
\begin{iii}
\item
The proof given above, shows that, in addition to the specialisations of the invertible sheaves
from~$\Pic S_{\overline{K}}$,
the Picard group
of~$\smash{\Pic S_{\overline\bbF_{\!\frakp}}}$
has (at least) two further~generators. One~of them may be chosen to be defined over
$\smash{\bbF_{\!\frakp}}$,
the other over its quadratic extension.
\item
Without the hypothesis on the determinant of the Frobenius, the argument simply reproves the standard fact that
$\rk \Pic S_{\overline\bbF_{\!\frakp}} \geq \rk \Pic S_{\overline{K}}$.
\end{iii}
\end{remso}

\begin{theoo}
\label{K3_fl}
Let\/~$K$
be a number field and\/
$S$
a\/
$K3$~surface
over\/~$K$.
Moreover,~let\/
$\frakp \subset \calO_K$
be a prime of residue characteristic\/
$\neq \!2$
and good reduction.

\begin{abc}
\item
Then~the following two equations hold,
\begin{align*}
\det(\Frob_\frakp \colon H^2_\et(S_{\overline{K}}, \bbQ_l(1)) \righttoleftarrow) &\!=\! \left( \frac{\Delta_{H^2}(S)}\frakp \right) \;\;\hspace{2cm} {\text and}\;\; \\ \det(\Frob_\frakp \colon T \righttoleftarrow) &\!=\! \left( \frac{\Delta_{H^2}(S)\Delta_{\Pic}(S)}\frakp \right) .
\end{align*}
\item
If\/
$\rk \Pic S_{\overline{K}}$
is even,~then
$$\left(\frac{\Delta_{H^2}(S)\Delta_{\Pic}(S)}\frakp\right) = -1 \quad\Longrightarrow\quad \rk \Pic S_{\overline\bbF_{\!\frakp}} \geq \rk \Pic S_{\overline{K}} + 2 \, .$$
In~other words, if\/
$\Delta_{H^2}(S)\Delta_{\Pic}(S)$
is not a square in\/
$K$~then
$$\{\,\frakp \mid \frakp \text{ \rm inert in } K(\sqrt{\Delta_{H^2}(S)\Delta_{\Pic}(S)})\,\} \subseteq \Pi_\jump(S) \, .$$
\end{abc}

\noindent
{\bf Proof.}
{\em
a)
The first formula is a particular case of formula~(\ref{Hi_char}). The~second one is a consequence of the first together with formulae (\ref{det_NS}) and~(\ref{det_Formel}).\smallskip

\noindent
b)
This follows from~a), together with Proposition~\ref{jumps}.
}
\eop
\end{theoo}

\begin{coroo}
\label{liminf}
Let\/~$K$
be a number field and\/
$S$
a\/
$K3$~surface
over\/~$K$.
Assume~that\/
$\Delta_{H^2}(S)\Delta_{\Pic}(S)$
is a non-square
in\/~$K$.
Then
$$\textstyle \liminf\limits_{B\to\infty} \gamma(S,B) \geq \frac12 \,.$$
\end{coroo}

\begin{defio}
For\/~$K$
a number field
and\/~$S$
a\/
\mbox{$K3$~surface}
over\/~$K$,
we call the quadratic~character\/
$$\smash{\tau_S := [\det H^2(S_{\overline{K}}, \bbQ_l(2))] \!\cdot\! [\det \Pic(S_{\overline{K}}) \otimes_\bbZ \bbQ] \colon \Gal(\overline{K}/K) \longrightarrow \{1,-1\}}$$
the {\em jump character\/}
of\/~$S$.\smallskip
\end{defio}

\begin{remso}
\begin{iii}
\item
The jump character
$\tau_S$
is given by
$$\frakp \mapsto \bigg(\frac{\Delta_{H^2}(S)\Delta_{\Pic}(S)}\frakp\bigg) \,,$$
for all good
primes~$\frakp$.
\item
Proposition~\ref{jumps} shows that, for
$S$
a
$K3$~surface
of even geometric Picard rank,
$\tau_S(\frakp) = -1$
implies
$\rk \Pic S_{\overline\bbF_{\!\frakp}} \geq \rk \Pic S_{\overline{K}} + 2$.
\end{iii}
\end{remso}

In this section, the assumption on the surface to be of type
$K3$
was used only in referring to the Tate conjecture. We~actually showed the~following.

\begin{theoo}
Let\/~$K$
be a number field and\/
$S$
a smooth and proper surface
over\/~$K$,
for which the Tate conjecture~holds. Moreover,~let\/
$\frakp \subset \calO_K$
be a prime of good reduction and suppose that the Tate conjecture holds
for~$S_{\bbF_{\!\frakp}}$,~too.
Then,~in the situation that\/
$\rk\NS S_{\overline{K}} \equiv \dim H^2_\et(S_{\overline{K}}, \bbQ_l(1)) \pmod 2$,
one~has
\begin{align*}
[\det H^2(S_{\overline{K}}, \bbQ_l(2))] \!\cdot\! [\det \Pic(S_{\overline{K}}) \otimes_\bbZ \bbQ](\frakp) = -1 \Longleftrightarrow& \\
\left(\frac{\Delta_{H^2}(S)\Delta_{\NS}(S)}\frakp\right) = -1 \Longrightarrow& \rk \NS S_{\overline\bbF_{\!\frakp}} \geq \rk \NS S_{\overline{K}} + 2 \, .
\end{align*}
\eop
\end{theoo}

\subsection{The criterion for non-triviality}

\begin{lemo}
\label{orth_sum}
Let\/~$p\neq2$
be a prime number and\/
$M$
a free\/
\mbox{$\bbZ_p$-module}
of finite
rank\/~$r$,
equipped with a non-degenerate symmetric~pairing. Let~an orthogonal map\/
$V\colon M \to M$
be given with characteristic polynomial\/
$(T-1)^{r_1}(T+1)^{r_2}$.
Then\/~$M$
is the orthogonal direct sum of the two generalised eigenspaces,
$$M = \ker (1-V)^{r_1} \oplus \ker (1+V)^{r_2} \,.$$
{\bf Proof.}
{\em
It is well known that the generalised eigenspaces for two distinct eigenvalues only have the zero element in~common. Moreover,
$\smash{\frac{1+V}2 + \frac{1-V}2= \id}$
yields
$$\sum\limits_{k=0}^r \textstyle \genfrac(){0pt}{1}{r}{k} (\frac{1+V}2)^k (\frac{1-V}2)^{r-k}(x) = x\,,$$
for
every~$x \in M$.
According to Cayley--Hamilton, the summands for
$k \geq r_2$
are contained in 
$\ker (1-V)^{r_1}$,
while those for
$k \leq r_2$,
i.e.\
$r-k \geq r_1$,
are contained in 
$\ker (1+V)^{r_2}$.
Thus,~the sum is the whole
of~$M$.
Finally,~it is a classical result for orthogonal maps that the generalised eigenspaces for two eigenvalues
$\lambda_1, \lambda_2$
with
$\lambda_1 \lambda_2 \neq 1$
are perpendicular. An~argument is given, e.g., in~\cite[Proposition~10.4.1]{Wal}.
}
\eop
\end{lemo}

\begin{theoo}[The vanishing cycle]
\label{K3_vancyc}
Let\/~$K$
be a number field
and\/~$S$
a\/
$K3$~surface
over\/~$K$.
Moreover,~let\/
$\frakp \subset \calO_K$
be a prime of residue characteristic\/
$\neq \!2$
such that\/
$S$~has
a regular, projective model\/
$\underline{S}$
over\/~$\calO_{K,\frakp}$,
the geometric fibre
$\underline{S}_{\overline\frakp}$
of which has exactly one singular
point\/~$\frakz$.
Assume\/~$\frakz$
to be an ordinary double~point.\smallskip

\noindent
Then,~for every prime
number\/~$l$,
different from the residue characteristic
of\/~$\frakp$,
the vanishing cycle\/~\cite[Exp.~XV, Th\'eor\`eme~3.4.(i)]{SGA7} associated
with~$\frakz$,
fulfils
$$\delta_{\frakz,l} \in H_{\alg,l}^\perp \, .$$
{\bf Proof.}
{\em
{\em First step.}
Generalities.

\noindent
Let~us denote the residue characteristic
of~$\frakp$
by~$p$.
On~the scheme
$S_{\overline{K}}$,
there is a monodromy automorphism \cite[Exp.~XV, Proposition~3.2.1.(ii)]{SGA7}, which is induced by a particular non-trivial element
$\smash{\nu \in I_\frakp \subset \Gal(\overline{K}/K)}$
of the inertia group. For~every prime
number~$l$,
including
$l=p$,
the induced map on
\mbox{$l$-adic}
cohomology is called the monodromy operator
$V \colon H^2_\et(S_{\overline{K}}, \bbZ_l(1)) \righttoleftarrow$.
This~is an orthogonal map with respect to the cup product~pairing. By~a slight abuse of notation, we denote the map induced
by~$\nu$
on
$\Pic S_{\overline{K}}$
by
$V$,~too.

If~$l\neq p$
then the action
of~$V$
is described by the Picard--Lefschetz~formula~\cite[Exp.~XV, Th\'e\-o\-r\`eme~3.4.(iii)]{SGA7}
\begin{align}
\label{PL}
V(c) = c + \langle c, \delta_{\frakz,l}\rangle\delta_{\frakz,l} \, .
\end{align}
Note~here that
$V$,
being induced by an element from the inertia group, acts trivially on
$\bbZ_l(1)$
itself. The~class
$\delta_{\frakz,l} \in H^2_\et(S_{\overline{K}}, \bbZ_l(1))$
is the so-called vanishing~cycle. It~is known that
$\langle \delta_{\frakz,l}, \delta_{\frakz,l}\rangle = -2$
\cite[Exp.~XV, Th\'e\-o\-r\`eme~3.4.(i)]{SGA7}.
In~particular, the Picard--Lefschetz~formula shows
$V(\delta_{\frakz,l}) = -\delta_{\frakz,l}$.
Moreover,~the operator
$V$
acts with characteristic polynomial
$(T-1)^{21}(T+1)$.

When~the action on
$H^2_\et(S_{\overline{K}}, \bbZ_p(1))$
is concerned, the characteristic polynomial is the same~\cite[Theorem~3.1]{Oc}, cf.\ \cite[\S\S2.3 and~2.4]{Fo}. It~seems, however, not to be known whether the action
of~$V$
is~semisimple. In~particular, no Picard--Lefschetz~formula is~available.\smallskip

\noindent
{\em Second step.}
If, for
some~$l_0\neq p$,
$\smash{\delta_{\frakz,l_0} \not\in H_{\alg,l_0}^\perp}$
then
$V$
acts non-trivially
on~$\Pic S_{\overline{K}}$.

\noindent
The~assumption means that
$\langle c, \delta_{\frakz,l_0}\rangle \neq 0$,
for a certain
class~$c \in H_{\alg,l_0}$.
Then~formula~(\ref{PL}) immediately implies that
$$\delta_{\frakz,l_0} = \frac1{\langle c, \delta_{\frakz,l_0}\rangle}[V(c)-c] \in H_{\alg,l_0} \, .$$
I.e.,~the operator
$V$
acts on
$H_{\alg,l_0} \cong \Pic S_{\overline{K}} \otimes_\bbZ \bbQ_{l_0}$
non-trivially with one eigenvector
$(-1)$,
while all others are equal
to~$1$.
As~$\pm1 \in \bbQ$,
the eigenspaces are defined in
$\Pic S_{\overline{K}} \otimes_\bbZ \bbQ$
already, which implies the~claim.%
%
\smallskip

\noindent
{\em Third step.}
Let
$d \in \Pic(S_{\overline{K}}) \subset \Pic S_{\overline{K}} \otimes_\bbZ \bbQ$
be a generator of the
\mbox{$(-1)$-eigenspace},
which is minimal, i.e.\ not divisible by any integer
$\neq\! \pm1$.
Then~$\langle d, d\rangle = -2$.

\noindent
For every prime
number~$l$,
the inclusion
$\smash{\Pic S_{\overline{K}}/l \stackrel{\overline{c}_1}{\hookrightarrow} H^2_\et(S_{\overline{K}}, \bbZ_l(1))/l}$
coming from the Kummer sequence shows that
$c_1(d) \in H^2_\et(S_{\overline{K}}, \bbZ_l(1))$
is not divisible
by~$l$.
Hence,~$c_1(d)$
is actually a generator of the
\mbox{$(-1)$-eigenspace}
in~$H^2_\et(S_{\overline{K}}, \bbZ_l(1))$.\smallskip

\noindent
{\em Claim:}
$\langle d, d\rangle$
is an
\mbox{$l$-adic}
unit, for every
prime~$l \neq 2$.

\noindent
Indeed,~for
$l \neq 2,p$,
the
\mbox{$\bbZ_p$-module}
$H^2_\et(S_{\overline{K}}, \bbZ_l(1))$
is the orthogonal direct sum of the
\mbox{$1$-}
and
\mbox{$(-1)$-eigenspaces},
since
$\smash{\frac12 \in \bbZ_l}$
and every cohomology
class~$c$
may be written in the form
$\frac12(c+V(c)) + \frac12(c-V(c))$.
For~$l=p$,
the
\mbox{$(-1)$-eigenspace}
is a direct summand, too, due to Lemma~\ref{orth_sum}.
As~the pairing on the total space is perfect, the same is true for the direct summands, which implies the~claim.\smallskip

\noindent
Thus,
$\langle d, d\rangle = \pm2^k$
for some non-negative
integer~$k$.
On~the other hand,
$c_1(d) \in H^2_\et(S_{\overline{K}}, \bbZ_2(1))$
is a generator of the
\mbox{$(-1)$-eigenspace}
and
$\delta_{\frakz,2}$
is another. Hence,
$c_1(d) = u \cdot \delta_{\frakz,2}$,
for a certain unit
$u \in \bbZ_2^*$.
Consequently,
$$\langle d, d\rangle = \langle c_1(d), c_1(d) \rangle = u^2 \langle \delta_{\frakz,2}, \delta_{\frakz,2} \rangle = -2u^2 \,,$$
which immediately shows that
$k=1$.
Moreover,~the minus sign is correct, since
$(-1)$
is not a square
in~$\bbQ_2$.\smallskip

\noindent
{\em Fourth step.}
Conclusion.

\noindent
In~the particular case of a
$K3$~surface,
it is well known that, for a class
$d \in \Pic(S_{\overline{K}})$
with
$\langle d, d\rangle = -2$,
either
$d$
or~$(-d)$
is represented by an effective divisor~\cite[Chap.~VIII, Proposition~3.7]{BHPV}. But~then
$V$
cannot interchange the two, a~contradiction.%
}%
\eop
\end{theoo}

\begin{remso}
\begin{iii}
\item
The regularity of the model
$\underline{S}$
implies that the singular point
on~$\underline{S}_{\overline\frakp}$
does not lift to a
$\calO_{K_\frakp}^\nr$-rational
point
on~$\underline{S}$.
\item
When there are two singular points instead of one, then the argument above then only shows that a non-trivial linear combination of
$\delta_{\frakz_1}$
and~$\delta_{\frakz_2}$
lies
in~$H_\alg^\perp$.
The~splitting field of
$\Pic(S_{\overline{K}}) \otimes_\bbZ \bbQ$
may well ramify then. Cf.~Corollary~\ref{Pic_split}.a),~below.
\item
There does not seem to be an obvious generalisation to other types of surfaces. For example, for rational surfaces one has
$\smash{H_\alg^\perp = 0}$,
but
$\smash{\delta_{\frakz,l}}$
is clearly nonzero. Also, the argument heavily relies on the fact that, for
$d \in \Pic(S_{\overline{K}})$
with
$\langle d, d\rangle = -2$,
either
$d$
or~$(-d)$
is effective, which seems to be rather specific for
$K3$
surfaces.
\end{iii}
\end{remso}

\begin{coroo}
\label{Pic_split}
Let\/~$K$
be a number field
and\/~$S$
a\/
$K3$~surface
over\/~$K$.
Moreover,~let\/
$\frakp \subset \calO_K$
be a prime of residue characteristic\/
$\neq \!2$
such that\/
$S$~has
a regular, projective model\/
$\underline{S}$
over\/~$\calO_{K,\frakp}$,
the geometric fibre
$\underline{S}_{\overline\frakp}$
of which has exactly one singular
point\/~$\frakz$.
Assume\/~$\frakz$
to be an ordinary double~point.

\begin{abc}
\item
Then~the splitting field of\/
$\Pic(S_{\overline{K}}) \otimes_\bbZ \bbQ$
is unramified
at\/~$\frakp$.
\item
The~character\/
$\smash{\big(\frac{\Delta_{\NS}(S)}\cdot\big) = [\det \NS(S_{\overline{K}}) \otimes_\bbZ \bbQ]}$
is unramified
at\/~$\frakp$.
Equivalently, the splitting field\/
$L_{S,\alg}$
is unramified
at\/~$\frakp$.
\end{abc}\smallskip

\noindent
{\bf Proof.}
{\em
a)
Choose a prime
number\/~$l$,
different from the residue characteristic
of\/~$\frakp$.
Then~there is the short exact sequence~\cite[Exp.~XV, Th\'eor\`eme~3.4.(ii)]{SGA7}
\begin{align*}
0 \longrightarrow  H^2_\et(\underline{S}_{\overline\frakp}, \bbQ_l(1)) \longrightarrow H^2_\et(S_{\overline{K}}, \bbQ_l(1)) &\longrightarrow \bbQ_l \longrightarrow 0 \,, \\[-1mm]
c \;&\mapsto\, \langle c, \delta_{\frakz,l}\rangle \,,
\end{align*}
provided by the theory of vanishing~cycles. Together~with the result of Theorem \ref{K3_vancyc}, it shows that every invertible sheaf on
$S_{\overline{K}}$
extends to 
$\underline{S}_{\overline\frakp}$.
In~other words, the splitting field
of~$\Pic(S_{\overline{K}}) \otimes_\bbZ \bbQ$
is contained
in~$K_\frakp^\nr$.\smallskip

\noindent
b)
This~is a direct consequence of~a).
}%
\eop
\end{coroo}

\begin{propo}[Reduction to one ordinary double point]
\label{red_A1}
Let\/~$K$
be a number field
and\/~$S$
a proper\/
\mbox{$K$-scheme}
that is pure of even
dimension\/~$i$.
Moreover,~let\/
$\frakp \subset \calO_K$
be a prime of residue characteristic\/
$\neq \!2$
such that\/
$S$~has
a regular, projective model\/
$\underline{S}$
over\/~$\calO_{K,\frakp}$,
the geometric fibre
$\underline{S}_{\overline\frakp}$
of which has exactly one singular
point\/~$\frakz$.
Assume\/~$\frakz$
to be an ordinary double~point.

\begin{abc}
\item
Then,~for any
prime\/~$l$
different from the residue characteristic
of\/~$\frakp$,
the\/
\mbox{$\smash{\Gal(\overline{K}/K)}$-representation\/}
$\smash{H^i_\et(S_{\overline{K}}, \bbQ_l(i/2))}$
is tamely ramified
at\/~$\frakp$.
The\/
\mbox{$\frakp$-adic}
valuation of its conductor is equal
to~$1$.
\item
The~quadratic character\/
$\smash{\big(\frac{\Delta_{H^i}(S)}\cdot\big) = [\det H^i(S_{\overline{K}}, \bbQ_l(i/2))]}$
is ramified
at\/~$\frakp$.
Equivalently, the splitting field\/
$L_S$
is ramified
at\/~$\frakp$.
\end{abc}\smallskip

\noindent
{\bf Proof.}
{\em
a)
In~this generality, the short exact sequence provided by the theory of vanishing cycles reads~\cite[Exp.~XV, Th\'eor\`eme~3.4.(ii)]{SGA7}
\begin{align*}
0 \longrightarrow  H^i_\et(\underline{S}_{\overline\frakp}, \bbQ_l(i/2)) \longrightarrow H^i_\et(S_{\overline{K}}, \bbQ_l(i/2)) &\longrightarrow \bbQ_l \longrightarrow 0 \,, \\[-1mm]
c \;&\mapsto\, \langle c, \delta_{\frakz,l}\rangle \,,
\end{align*}
with~$\smash{(\delta_{\frakz,l}, \delta_{\frakz,l}) = (-1)^{i/2} \!\cdot\! 2}$.
The
\mbox{$\frakp$-adic}
valuation of the conductor is determined by the restriction of the representation to the decomposition group
$\smash{D_\frakp \cong \Gal(\overline{K}_\frakp/K_\frakp)}$.
The~exact sequence shows that the subspace
$\smash{\delta_{\frakz,l}^\perp \subset H^i_\et(S_{\overline{K}}, \bbQ_l(i/2))}$
is acted upon via the quotient
$\smash{\Gal(\overline\bbF_{\!\frakp}/\bbF_{\!\frakp}) \cong D_\frakp/I_\frakp}$.
I.e.,~the inertia
group~$I_\frakp$
fixes
$\smash{\delta_{\frakz,l}^\perp}$
pointwise.
As~one has
$V(\delta_{\frakz,l}) = -\delta_{\frakz,l}$
for
$V \in I_\frakp$,
due to the Picard--Lefschetz~formula, this~yields
$$H^i_\et(S_{\overline{K}}, \bbQ_l(i/2))^{I_\frakp} = \delta_{\frakz,l}^\perp \,.$$
Moreover,~the action
of~$I_\frakp$
respects orthogonality and cup product pairing, so
$\delta_{\frakz,l}$
can be mapped only
to~$\pm \delta_{\frakz,l}$.
Thus,~there is a subgroup
$I' \subset I_\frakp$
of index two acting~trivially. Since~the residue characteristic
of~$\frakp$
is~$\neq \!2$,
this yields~tameness.

In~this case, the
\mbox{$\frakp$-adic}
valuation of the conductor is given by~\cite[formulae~(11) and~(8)]{Se}
\begin{align*}
\dim_{\bbQ_l} H^i_\et(S_{\overline{K}}, \bbQ_l(i/2)) - \dim_{\bbQ_l} H^i_\et(S_{\overline{K}}, \bbQ_l(i/2))^{I_\frakp} &= \\[-1.5mm]
\dim_{\bbQ_l} H^i_\et(S_{\overline{K}}&, \bbQ_l(i/2)) - \dim_{\bbQ_l} \delta_{\frakz,l}^\perp = 1 \,.
\end{align*}

\noindent
b)
This~is very easily shown~directly. The~monodromy operator
$V \in I_\frakp$
fixes all cohomology classes perpendicular
to~$\delta_{\frakz,l}$
and sends
$\delta_{\frakz,l}$
to~$(-\delta_{\frakz,l})$.
Therefore,
$\det (V\colon H^2_\et(S_{\overline{K}}, \bbQ_l(1)) \righttoleftarrow) = -1$.
In~particular, 
$\ker\, [\det H^2(S_{\overline{K}}, \bbQ_l(1))]$
does not include all
of~$I_\frakp$,
and hence the field corresponding under the Galois correspondence to
$\ker\, [\det H^2(S_{\overline{K}}, \bbQ_l(1))]$
is not contained in
$K_\frakp^\nr$.
}
\eop
\end{propo}

\begin{remo}
Suppose that, for some prime
$\frakp$,
$S$
has a model of bad reduction of the kind as described~above.
Then~$S$
does not have a model of good
reduction~at~$\frakp$.
Indeed,~the conclusions of Lemma~\ref{unverzw_Hi}.b) and Proposition~\ref{red_A1}.b) are independent of the model. The~existence of a model of good reduction implies that
$\smash{\big(\frac{\Delta_{H^i}(S)}\cdot\big)}$
is unramified, while the existence of a model of the type above enforces~ramification.
\end{remo}

\begin{coroo}[The jump character]
\label{jchar_nontriv}
Let\/~$K$
be a number field
and\/~$S$
a\/
$K3$~surface
over\/~$K$.
Moreover,~let\/
$\frakp \subset \calO_K$
be a prime of residue characteristic different
from\/~$2$.

\begin{abc}
\item
If\/~$S$
has good reduction
at\/~$\frakp$,
then\/
$\smash{\tau_S = \big(\frac{\Delta_{H^2}(S)\Delta_{\Pic}(S)}\cdot\big)}$
is unramified
at\/~$\frakp$.
\item
Suppose~that\/
$S$
has a regular, projective model\/
$\underline{S}$
over\/~$\calO_{K,\frakp}$,
the geometric fibre
$\underline{S}_{\overline\frakp}$
of which has exactly one singular point, and assume this to be an ordinary double~point. Then~the jump~character\/
$\smash{\tau_S = \big(\frac{\Delta_{H^2}(S)\Delta_{\Pic}(S)}\cdot\big)}$
ramifies
at\/~$\frakp$.
\end{abc}

\noindent
{\bf Proof}.
{\em
a)~is clear from Lemmata~\ref{unverzw_Hi}.b) and~\ref{unverzw}.b), while the assertion of part~b) follows from Corollary~\ref{Pic_split}.b) together with Proposition~\ref{red_A1}.b).
}
\eop
\end{coroo}

\subsection{Examples and experimental results}\leavevmode

\begin{algoo}[Computing~$\Delta_{H^2}(S)$]
\label{alg_delta}
Given~a proper surface
$S$
over~$\bbQ$,
the set
$\{q_1, \ldots, q_m\}$
of all bad primes
of~$S$,
and an oracle for
$\smash{\det(\Frob_p\colon H^2_\et(S_{\overline\bbQ}, \bbQ_l(1)) \righttoleftarrow)}$
for any
$p \neq q_j$,
this algorithm computes
$\Delta_{H^2}(S)$.

\begin{iii}
\item
Add
$q_0 := -1$
to the list of bad primes.
\item
Build a matrix
$A$,
the entries of which are the Legendre symbols
$\smash{(\frac{q_j}{p_i})}$,
for good primes
$p_i$
chosen at random. Keep adding rows until the matrix has
rank~$m+1$
over
$\{1,-1\} \cong \bbZ/2\bbZ$.
\item
Put
$\smash{b_i = \det(\Frob_{p_i}\colon H^2_\et(S_{\overline\bbQ}, \bbQ_l(1)) \righttoleftarrow)}$
and solve the linear system
$Ax=b$
of equations. If~the solution vector is
$(x_0, \ldots, x_m) \in (\bbZ/2\bbZ)^m$
then
$\Delta_{H^2}(S)$
is the class of
$(-1)^{e_0} q_1^{e_1}\cdots q_m^{e_m}$
in~$\bbQ^*/(\bbQ^*)^2$,
for
$e_i \in \{0,1\} \subset \bbZ$
representing the residue class
$x_i \in \bbZ/2\bbZ$.
\end{iii}
\end{algoo}

\begin{remso}
\begin{iii}
\item\looseness-1
The oracle for
$\smash{\det(\Frob_p\colon H^2_\et(S_{\overline\bbQ}, \bbQ_l(1)) \righttoleftarrow)}$
is, of course, provided by counting the points
on~$S$
that are defined over
$\bbF_{\!p}$
and some of its~extensions.
\item
Dirichlet's Theorem on primes in arithmetic progressions ensures that there exist primes so that the matrix
$A$
has
rank~$m+1$.
\item
(An improvement.)
In the case that
$2$
is a good prime, step~i) of Algorithm~\ref{alg_delta} may be omitted. The~sign
of~$\Delta_{H^2}(S)$
is then determined by the condition that the character be unramified
at~$2$.
\item
(A further improvement.)
Assume~that the 
surface~$S$
is~$K3$.
Then,~for some or many of its bad primes
$p \neq 2$,
it may happen that Proposition~\ref{red_A1}.b)~applies. At~such a prime, the jump character necessarily ramifies, which means that
$\Delta_{H^2}(S)$
must be of odd
\mbox{$p$-adic}
valuation. Thus,~the solution vector is bound to have a
component~$1 \in \bbZ/2\bbZ$
at the corresponding~coordinates.

If~$S$
has a model in some
$\Pb^N_\bbZ$
that is given by explicit equations, one may compute the set of all bad primes using Gr\"obner bases and integer factorisation, and finally analyse the singular~points. 
Having found several primes to which Proposition~\ref{red_A1}.b) applies, this information may be used in order to get by with a linear system of equations of smaller~size. In~other words, less point counting is~necessary.

It~is our experience that this improvement of Algorithm~\ref{alg_delta} often leads to an enormous gain for a ``random'' surface, while for the constructed examples, which we present below, it would not help~much.
\item
There is an obvious modification of Algorithm~\ref{alg_delta} to directly determine the jump character.
\end{iii}
\end{remso}

\begin{algoo}[Statistical algorithm computing the jump character]
\label{alg_jump}
Given a
$K3$~surface
$S$
over~$\bbQ$
of even geometric Picard~rank, the set
$\{q_1, \ldots, q_m\}$
of all bad primes
of~$S$,
and a list of good non-jump primes, this algorithm determines a finite subgroup, containing the jump character, of the group of all characters of
$\Gal(\overline\bbQ/\bbQ)$
with values
in~$\{1, -1\}$.

\begin{iii}
\item
Add
$q_0 := -1$
to the list of bad primes.
\item
Build a matrix
$A$,
the entries of which are the Legendre symbols
$\smash{(\frac{q_j}{p_i})}$,
for the non-jump
primes~$p_i$.
\item
Determine~the kernel
of~$A$.
From~each kernel vector, calculate a candidate for the jump character in the same way as in Algorithm~\ref{alg_delta}.iii).
\end{iii}
\end{algoo}

\begin{remso}
\begin{iii}
\item
If~the kernel is the zero space, then this proves the jump character to be trivial. If~the kernel is one-dimensional, then there are two possible answers. A non-trivial character, which is directly computed from a kernel vector, and the trivial~one.
\item
If~the kernel is still one-dimensional when the system of equations is rather overdetermined, then this gives strong evidence for the jump character to be non-trivial. In~practice, we work with at least
$4(m+1)$
non-jump~primes.
\item
The~trivial character is unramified at every~prime. Thus,~as~soon as it applies, Corollary~\ref{jchar_nontriv}.b) excludes the trivial character, and therefore makes the outcome of Algorithm~\ref{alg_jump} usually~unique. Corollary~\ref{jchar_nontriv}.b)~is useful as well to accelerate the~calculations.
\end{iii}
\end{remso}

\begin{exo}
\label{deg4}
Let~$S$
be the diagonal quartic
in~$\Pb^3_\bbQ$,
given by
$X_0^4 + X_1^4 + X_2^4 + X_3^4 = 0$.
Then the geometric Picard rank
of~$S$
is~$20$
and the jump character is given by
$\smash{(\frac{-1}{\cdot})}$.\smallskip

\noindent
{\bf Proof.}
The model
$\calS$
of
$S$
that is given in
$\Pb^3_\bbZ$
by the same equation has bad reduction only
at~$2$.
Hence,~$\Delta_{H^2}(S) = \pm1$
or~$\pm2$.
Counting points on the reductions
$\smash{S_{\bbF_{\!3}}}$
and~$\smash{S_{\bbF_{\!5}}}$,
one finds that
$\smash{\det(\Frob_p \colon H^2_\et(S_{\overline\bbQ}, \bbQ_l(1)) \righttoleftarrow) = 1}$
for both
$p=3$
and~$5$.
Thus,~Algorithm \ref{alg_delta} shows that
$\smash{\Delta_{H^2}(S) = 1}$.

On~the other hand, it is classically known that the 48 lines
on~$\smash{S_{\overline\bbQ}}$
generate the geometric Picard group, which is of 
rank~$20$.
In~particular,
$\smash{\Pic S_{\overline\bbQ}}$
is defined over
$\bbQ(\zeta_8) = \bbQ(i, \sqrt{2})$.
Moreover,~\cite[Appendix A, Examples A62, B33, C27, and D27]{Br} show that the Galois representation
$\Pic(S_{\overline\bbQ}) \otimes_\bbZ \bbC$
splits into characters as
$$\chi_\triv^5 \oplus \chi_{\bbQ(i)}^3 \oplus \chi_{\bbQ(\sqrt{2})}^6 \oplus \chi_{\bbQ(\sqrt{-2})}^6 \, .$$
Here,~for
$K$
a quadratic number field,
$\chi_K\colon \Gal(\overline\bbQ/\bbQ) \to \{1,-1\}$
denotes the character that becomes trivial after restriction to
$\Gal(\overline\bbQ/K)$
and defines the non-trivial quadratic character on
$\Gal(K/\bbQ)$.
Consequently,
$$\Exterior^\maxi \Pic(S) \otimes_\bbZ \bbC = \chi_\triv^{\otimes5} \otimes \chi_{\bbQ(i)}^{\otimes3} \otimes \chi_{\bbQ(\sqrt{2})}^{\otimes6} \otimes \chi_{\bbQ(\sqrt{-2})}^{\otimes6} = \chi_{\bbQ(i)} \,,$$
$\Delta_{\Pic}(S) = -1$
and~$\Delta_{H^2}(S)\Delta_{\Pic}(S) = -1$.
\eop
\end{exo}

\begin{remo}
It is known at least since 1963 \cite{T} that, in this example, there are no rank jumps, except for those explained by the jump character. I.e.,~one has
$\smash{\rk \Pic S_{\overline\bbF_{\!p}} = 20}$
for all primes
$p \equiv 1 \pmod 4$.
In~fact, the eigenvalues of
$\Frob_p$
on
$\smash{H^2_\et(S_{\overline\bbQ}, \bbQ_l(1))}$
may be determined using Jacobi sums \cite[Chapter~8, Theorem~5]{IR} and it turns out that two of them are
$\smash{\frac{\pi^2}p}$
and its conjugate, for
$p = \pi\overline\pi$
a factorisation
in~$\bbQ(i)$.
Cf.~\cite[particularly formulae (12) and~(13)]{PS} for more~details.
\end{remo}

\begin{exo}
\label{deg2}
Let~$S$
be the double cover
of~$\Pb^2_\bbQ$,
given by
$w^2 = X_0^6 + X_1^6 + X_2^6$.
Then the geometric Picard rank
of~$S$
is~$20$
and the jump character is given by
$\smash{(\frac{-3}{\cdot})}$.\smallskip

\noindent
{\bf Proof.}
The double cover
$\calS$
of~$\Pb^2_\bbZ$
that is given by the same equation has bad reduction only at the primes
$2$
and~$3$.
Hence,~$\smash{\Delta_{H^2}(S) = \pm1}$,
$\smash{\pm2}$,
$\smash{\pm3}$
or~$\smash{\pm6}$.
Furthermore,~counting points on the reductions
$S_{\bbF_{\!5}}$,
$S_{\bbF_{\!7}}$,
and~$S_{\bbF_{\!13}}$,
one finds that
$\det(\Frob_p \colon H^2_\et(S_{\overline\bbQ}, \bbQ_l(1)) \righttoleftarrow) = 1$,
for both
$p=5$
and
$p=13$,
and
$(-1)$
for
$p=7$.
Thus,~Algorithm \ref{alg_delta} shows
$\Delta_{H^2}(S) = -1$.

The ramification sextic allows 18 tritangent lines of the type
$X_i + \zeta_{12}^m X_j = 0$,
for
$m$~odd.
Furthermore,~the 18 conics of type
$\smash{X_i X_j + \frac{\zeta_6^m}{\sqrt[3]{2}} X_k^2 = 0}$
are six times tangent to the ramification~sextic.
The~irreducible components of the preimages of these curves together generate the geometric Picard group up to finite index. Working with the tritangents alone, one would end up with a sublattice that is not of full~rank.

We~implemented in {\tt magma} a function to compute intersection numbers
on~$S$
and, starting with 14 tritangent lines and six conics being six times tangent, found a non-degenerate
$20 \times 20$
intersection~matrix. Using this, it turns out that the splitting field of
$\smash{\Pic S_{\overline\bbQ}}$
is in fact
$\bbQ(\zeta_3, \sqrt[3]{2}, i)$,
having Galois group
$S_3 \times \bbZ/2\bbZ$.
Moreover,~the Galois representation
$\Pic(S_{\overline\bbQ}) \otimes_\bbZ \bbC$
splits into irreducible components as
$$\chi_\triv^4 \oplus \chi_{\bbQ(i)}^4 \oplus \chi_{\bbQ(\sqrt{3})}^3 \oplus \chi_{\bbQ(\sqrt{-3})}^3 \oplus V^3 \, ,$$
where the characters are defined as above and
$V$
denotes the irreducible two-dimensional representation of the factor group
$\smash{\Gal(\bbQ(\zeta_3, \sqrt[3]{2})/\bbQ) \cong S_3}$.
Consequently,
$$\Exterior^\maxi \Pic(S_{\overline\bbQ}) \otimes_\bbZ \bbC = \chi_\triv^{\otimes4} \otimes \chi_{\bbQ(i)}^{\otimes4} \otimes \chi_{\bbQ(\sqrt{3})}^{\otimes3} \otimes \chi_{\bbQ(\sqrt{-3})}^{\otimes3} \otimes \chi_{\bbQ(\sqrt{-3})}^{\otimes3} = \chi_{\bbQ(\sqrt{3})} \,,$$
which implies that
$\Delta_{\Pic}(S) = 3$
and~$\Delta_{H^2}(S)\Delta_{\Pic}(S) = -3$.
\eop
\end{exo}

\begin{remo}
Again, there are no rank jumps, except for those explained by the jump character. I.e.,~one has
$\smash{\rk \Pic S_{\overline\bbF_{\!p}} = 20}$
for all primes
$p \equiv 1 \pmod 3$.
The~eigenvalues of
$\Frob_p$
may again be determined using Jacobi~sums. Here,~it turns out that two of them are
$J(\omega,\omega,\omega)/p$
\cite[Proposition 8.5.1]{IR} and its conjugate, for
$\omega$~a
primitive sextic character
on~$\bbF_{\!p}^*$.
A~short calculation, using~\cite[Chapter~8, Theorem 3]{IR} and \cite[Theorems 3.1 and~3.4)]{BE} shows that these quantities evaluate to
$\smash{(\frac{-1}p) \frac{\pi^2}p}$
and its conjugate, for
$\pi$
a primary element
\cite[Proposition~9.3.5]{IR}
in~$\bbZ(\zeta_3)$
of
norm~$p$.
\end{remo}

\begin{remso}
\begin{iii}
\item
The surfaces described in Examples~\ref{deg4} and~\ref{deg2} are
$K3$
surfaces defined
over~$\bbQ$
of geometric Picard
rank~$20$,
and, as such, very particular~objects. Due~to the work of R.\ Livn\'e~\cite{Lv}, it is known that they are modular. Moreover,~there is a formula describing the determinant of Frobenius on the transcendental part of cohomology in terms of the discriminant of the Picard lattice alone \cite[Example~1.6]{Lv}. Our~calculations presented above are supposed to illustrate our method to compute the jump character in a situation where the Picard group is completely~known.
\item
The~surface from Example~\ref{deg2} is, up to isogeny, the Kummer surface associated with a product of an elliptic curve with itself.
(Over~$\bbC$,
this is classically known, cf.~\cite[Proof of Theorem~4]{SI}. 
Over~$\bbQ$,
an argument is given, e.g.,
in~\cite[Proposition~2.3]{BCFNW}.) Its~jump character may also be determined that~way.
\end{iii}
\end{remso}

\begin{exo}
\label{Kummer}
Let\/~$K$
be a number field and\/
$S$
the Kummer surface of an abelian surface
over\/~$K$
that geometrically splits into a product\/
$E_1 \!\times\! E_2$
of elliptic~curves. Assume~that\/
$\rk\Pic S_{\overline{K}} = 18$.
Then~there are two~cases.

\begin{abc}
\item
If~the elliptic curves\/
$E_1$
and\/~$E_2$
are defined
over\/~$K$
then the jump character
of\/~$S$
is~trivial.
\item
If~the elliptic curves\/
$E_1$
and\/~$E_2$
are defined over a quadratic extension\/
$K(\sqrt{d})$~and
conjugate to each other, then the jump character
of\/~$S$
is\/~$(\frac{d}.)$
\end{abc}\smallskip

\noindent
{\bf Proof.}
a)
The transcendental part
$T \subset H^2_\et(S_{\overline{K}}, \bbQ_l(1))$
is isomorphic~to
$$T \cong H^1_\et(E_1, \bbQ_l) \otimes H^1_\et(E_2, \bbQ_l(1)) \, ,$$
hence
\begin{eqnarray*}
\Exterior^\maxi T &\cong& \Exterior^{\!2} H^1_\et(E_1, \bbQ_l)^{\otimes 2} \otimes \Exterior^{\!2} H^1_\et(E_2, \bbQ_l(1))^{\otimes 2} \\
 &\cong& H^2_\et(E_1, \bbQ_l)^{\otimes 2} \otimes H^2_\et(E_2, \bbQ_l(2))^{\otimes 2} \\
 &\cong& H^2_\et(E_1, \bbQ_l(1))^{\otimes 2} \otimes H^2_\et(E_2, \bbQ_l(1))^{\otimes 2} \,,
\end{eqnarray*}
and both factors are acted upon trivially by
$\Gal(\overline{K}/K)$.\smallskip

\noindent
b)
Let
$\sigma \in \Gal(\overline{K}/K)$
be any automorphism that changes the sign
of~$\sqrt{d}$.
Then~$\sigma$
interchanges the components of
$H^1_\et(E_1 \!\times\! E_2, \bbQ_l) = H^1_\et(E_1, \bbQ_l) \oplus H^1_\et(E_2, \bbQ_l)$.
I.e.,
$\sigma$
acts with eigenvalues
$(-1)$
and~$1$,
both of
multiplicity~$2$.
Hence,~on
$H^2_\et(E_1 \!\times\! E_2, \bbQ_l) \cong \Exterior^{\!2} H^1_\et(E_1 \!\times\! E_2, \bbQ_l)$,
one has the eigenvalues
$(-1)$,
of
multiplicity~$4$,
and~$1$,
of
multiplicity~$2$.

However,~under
$\sigma$,
the two algebraic classes in
$H^2_\et(E_1, \bbQ_l) \oplus H^0_\et(E_2, \bbQ_l)$
and
$H^0_\et(E_1, \bbQ_l) \oplus H^2_\et(E_2, \bbQ_l)$
are interchanged, so that the eigenvalues
$(-1)$
and~$1$
occur on the algebraic part. Therefore,~the eigenvalues
on~$T$
are
$(-1)$,
with
multiplicity~$3$,
and
$1$,
with
multiplicity~$1$.
Hence,~every
$\sigma \in \Gal(\overline{K}/K)$
as chosen above acts as
$(-1)$
on~$\Exterior^\maxi T$,
which is enough to imply the~claim.
\eop
\end{exo}

\begin{exo}
In~\cite[Examples~3.3,~3.4,~and~3.5]{CT}, Yu.~Tschinkel and the first author reported numerical evidence for
$\liminf_{B \rightarrow \infty} \gamma(S, B) \geq 1/2$,
in the case of three
$K3$~surfaces
over
$\bbQ$
of geometric Picard rank~two.
    
This~indeed follows from Corollary~\ref{liminf}, once one proves that  $\smash{\Delta_{H^2}(S) \Delta_{\Pic}(S)}$ 
is not a square
in~$\bbQ$.
For~each of the examples, one has
$\smash{\Pic S_{\overline\bbQ} = \Pic S}$
and therefore
$\smash{\Delta_{\Pic}(S) = 1 \in \bbQ^*/(\bbQ^*)^2}$.
Moreover,~Algorithm~\ref{alg_delta} determines
$\Delta_{H^2}(S)$
to, in this~order,

~
{\vspace{-3mm}\tiny\interdisplaylinepenalty100
\begin{align*}
& -5 
\cdot 151 
\cdot 22490817357414371041 \cdot \\ 
&  387308497430149337233666358807996260780875056740850984213276970343278935342068889706146733313789 \, , \\[2mm]
& {} 53 
\cdot 2624174618795407 
\cdot 512854561846964817139494202072778341 \cdot \\
& 1215218370089028769076718102126921744353362873 \cdot \\
& 6847124397158950456921300435158115445627072734996149041990563857503 \, , \\[1mm]
& \text{and} \\[1mm]
& -47 
\cdot 3109 
\cdot 4969 
\cdot 14857095849982608071 
\cdot 445410277660928347762586764331874432202584688016149 \cdot \\
& 658652708525052699993424198738842485998115218667979560362214198830101650254490711 \, .
\end{align*}}%
Each of the factors listed is reported as being prime by~{\tt magma}, version 2.21.8.
\end{exo}

\begin{exo}
\label{jctriv}
Let\/~$S$
be the\/
$K3$~surface
over\/~$\bbQ$,
given by the equation
\begin{equation}
\label{spec_quart}
X_3^4 + f_2(X_0,X_1,X_2)X_3^2 + f_4(X_0,X_1,X_2) = 0 \, ,
\end{equation}
for
\begin{eqnarray*}
f_2(X_0,X_1,X_2) & := & X_0^2 - X_0X_1 - X_0X_2 - X_1X_2 \qquad {\rm and} \\
f_4(X_0,X_1,X_2) & := & - X_0^3X_2 + X_0X_1^2X_2 - X_1^4 - X_2^4 \, .
\end{eqnarray*}
Then~the geometric Picard rank
of\/~$S$
is\/~$8$
and the jump character
of~$S$
is~trivial.\medskip

\noindent
{\bf Proof.}
First~of all, a space
quartic~$S$
of the form~(\ref{spec_quart}) is of geometric Picard rank at
least~$8$.
Indeed,~the surface
$S$
comes equipped with a finite morphism
$p\colon S \to S'$,
which is generically
$2\!:\!1$,
to an underlying degree two del Pezzo
surface~$S'$.
The~induced homomorphism
$\smash{p^*\colon \Pic S'_{\overline{K}} \to \Pic S_{\overline{K}}}$
doubles all intersection numbers. As,~on a degree two del Pezzo surface, there are no non-trivial invertible sheaves that are numerically equivalent to zero, we see that
$p^*$
is necessarily injective. The~claim~follows.%

Thus,~for the first assertion, it suffices to find a
prime~$p$
of good reduction such that
$\smash{\rk \Pic S_{\overline\bbF_{\!p}} = 8}$.
For~example,
$p = 19$,
$43$,
$61$,
$101$,
$109$,
$139$,
$149$,
$151$,
$157$,
and~$163$
do the job, as is easily shown in the usual way, based on counting~points. Cf.~\cite{CT} for more details and further~references.

On~the other hand, a calculation using Gr\"obner bases shows that the model
$\calS$
of~$S$
given by the same equation
in~$\Pb^3_\bbZ$
has bad reduction only at the primes
$2$,
$3$,
$47$,
and~$431$.
Using~Algorithm~\ref{alg_jump}, one then proves the triviality of the jump~character. In~fact, only the first five non-jump primes
$19$,
$43$,
$61$,
$101$
and~$109$
are needed in order to do~this.
\eop
\end{exo}

\begin{remo}
This~example, and several others of the same kind, were found by a systematic inspection of all space quartics of the form~(\ref{spec_quart}), with coefficients from
$\{-1,0,1\}$.
This~led to a sample of
$183\,098\,318$
non-singular surfaces in total, among which only a few hundred have trivial jump character, together with geometric Picard
rank~$8$.
\end{remo}

\begin{exo}
\label{jcnontriv}
Let\/~$S$
be the\/
$K3$~surface
over\/~$\bbQ$,
given by the equation
\begin{equation}
\label{4_cover}
X_3^4 + f_4(X_0,X_1,X_2) = 0 \, ,
\end{equation}
for
$$f_4(X_0,X_1,X_2) := X_0^4 - X_0^3X_1 - 2X_0^3X_2 - X_0^2X_1X_2 + X_0X_1^2X_2 - X_1^4 - X_2^4 \, .$$
Then~the geometric Picard rank
of\/~$S$
is\/~$8$
and the jump character
of~$S$
is
$(\frac{-1}.)$.\medskip

\noindent
{\bf Proof.}
Again, for the first claim, it suffices to find a prime
$p$
of good reduction such that
$\smash{\rk \Pic S_{\overline\bbF_{\!p}} = 8}$.
For~example,
$p = 5$,
$13$,
$41$,
$53$,
$61$,
$73$,
$89$,
$97$,
$101$,
$109$,
$113$,
$137$,
$149$,
$157$,
$173$,
$181$,
$193$,
and
$197$
do the job.

Moreover, a calculation using Gr\"obner bases shows that the model
$\calS$
of~$S$
given
in~$\Pb^3_\bbZ$
by equation~(\ref{4_cover}) has bad reduction only at the primes
$2$,
$7$,
$6449$,
and~$39\,870\,353$.
For~the obvious integral model of the underlying degree two del Pezzo
surface~$S'$,
the same is~true.
Algorithm~\ref{alg_delta} then proves that
$\Delta_{H^2}(S) = 2 \cdot 7 \cdot 6449 \cdot 39\,870\,353$
and that
$\Delta_{H^2}(S') = -2 \cdot 7 \cdot 6449 \cdot 39\,870\,353$.
The~seven good primes up
to~$23$
are in fact~sufficient.

As, on a del Pezzo surface, every cohomology class is algebraic, we conclude that
$\Delta_{\NS}(S') = \Delta_{H^2}(S')$.
Furthermore,~the linear~map
$$p^* \!\otimes_\bbZ\! \bbQ\colon \NS(S') \!\otimes_\bbZ\! \bbQ = \Pic S' \!\otimes_\bbZ\! \bbQ \to \Pic S \!\otimes_\bbZ\! \bbQ$$
is an isomorphism, since it is injective and either
\mbox{$\bbQ$-vector}
space is of
dimension~$8$.
Therefore,~$\Delta_{\Pic}(S) = \Delta_{\NS}(S')$,
which implies the~claim.
\eop
\end{exo}

\begin{remo}
Note that, for the surface above, Algorithm~\ref{alg_jump} would only prove that the jump character
is
$(\frac{-1}.)$
or~trivial. Moreover,~the non-triviality criterion, given in Corollary~\ref{jchar_nontriv}.b), could not resolve the ambiguity~either.
%
\end{remo}

\subsection{Interaction of jumps}\leavevmode

\noindent
As is well known, the geometric Picard rank always jumps under reduction when
$\rk \Pic S_{\overline{K}}$
is odd. The same is true, when there is real multiplication by some
field~$E$
and one has an odd quotient
$(22 - \rk \Pic S_{\overline{K}})/[E\!:\!\bbQ]$
\cite[Theorem~1(2)]{Ch14}. One~might speculate in these cases, whether the jump character causes the jumps to be even~larger. This~does, however, not happen, as is shown by the examples~below.

\begin{lemo}[$K3$~surfaces
having a non-singular degree two model]
Let\/~$K$
be a number field and\/
$S$
a
$K3$~surface
over\/~$K$,
given by\/
$w^2 = f_6(X_0, X_1, X_2)$,
for\/
$f_6$
a homogeneous form of
degree\/~$6$.
Write\/
$$S_\lambda\colon \lambda w^2 = f_6(X_0, X_1, X_2)$$
for the quadratic twist by\/
$\lambda \in K^*$.
Then
$$\Delta_{H^2}(S_\lambda) = \lambda \Delta_{H^2}(S) \quad \text{and} \quad \Delta_{\Pic}(S_\lambda) = \lambda^{\rk\Pic S_{\overline{K}}-1} \Delta_{\Pic}(S) \, .$$

\noindent
{\bf Proof.}
{\em
Let
$\frakp$
be a good prime
of~$S$
such that
$\lambda$
is a
\mbox{$\frakp$-adic}
unit. Then,~for the reductions mod
$\frakp$,
one has that
$(S_\lambda)_\frakp$
is a non-trivial twist of
$S_\frakp$
in the case that
$\lambda$
is a non-square modulo
$\frakp$,
and
$S_\frakp \cong (S_\lambda)_\frakp$,
otherwise. The~assertion therefore follows from \cite[Fact~25]{EJ10}.
}
\eop
\end{lemo}

\begin{remo}[The odd rank case]
\label{comb_odd}
Assume that
$\rk\Pic S_{\overline{K}} = 1$.
Then,~for any
prime~$\frakp$
of good reduction, there exists some
\mbox{$\frakp$-adic}
unit
$\lambda \in K^*$
such that
$\Delta_{H^2}(S_\lambda) \Delta_{\Pic}(S_\lambda)$
is a non-square
modulo~$\frakp$.
If the effect of the odd rank added up with that of the transcendental character, then this would imply
$$\rk\Pic S_{\overline\bbF_{\!\frakp}} = \rk\Pic (S_\lambda)_{\overline\bbF_{\!\frakp}} \geq 4 \,.$$
There~are, however, explicit
degree~$2$
$K3$~surfaces
known of geometric Picard
rank~$1$
that reduce to geometric Picard
rank~$2$
at certain primes~\cite[Theorem~3.1]{vL}, cf.~\cite[Example~5.1.1]{EJ11}.
\end{remo}

\begin{exo}[The case of real multiplication]
\label{comb_rm}
Let\/~$S$
be the minimal desingularisation of the double cover
of~$\Pb^2_\bbQ$,
given by
$w^2 = X_0 X_1 X_2 \!\cdot\! f_3(X_0, X_1, X_2)$,
for
\begin{align*}
f_3(X_0, X_1, X_2) := X_0^3 + 3X_0^2X_1 - 2X_0^2X_2 + 5X_0X_1^2 - X_0X_2^2 + 3X_1^3 & \\[-1mm]
{} - 2X_1^2X_2 - & 3X_1X_2^2 + 2X_2^3 \, .
\end{align*}
There is strong evidence that
$S$
has real multiplication
by~$\bbQ(\sqrt{3})$.
Indeed,~$S$
is the surface
$\smash{V^{(3)}_{1,1}}$
from \cite[Conjectures~5.2]{EJ16}. Its~model
$\calS$
being the double cover
of~$\Pb^2_\bbZ$
given by the same equation has bad reduction only at
$2$,
$3$,
and~$5$.
Modulo~all other primes
$p < 1000$,
the reduction
$\calS_{\bbF_{\!p}}$
is of geometric Picard rank
$18$,
except for
$p = 263$,
where the geometric Picard rank
is~$22$.
On~the other hand, a sublattice of
$\smash{\Pic S_{\overline\bbQ}}$
of
rank~$16$
may be explicitly given. Altogether,~taking real multiplication for granted, one concludes that
$\smash{\rk\Pic S_{\overline\bbQ} = 16}$.

Concerning~$\smash{\Pic S_{\overline\bbQ}}$,
there are 13 obvious generators, given by the pull-back of a general line
on~$\smash{\Pb^2_\bbQ}$
and the exceptional curves obtained by blowing up the twelve singular points of the ramification~locus. Ten of these singular points are defined
over~$\bbQ$,
the other two
over~$\smash{\bbQ(\sqrt{-2})}$.
Hence,~this part of 
$\smash{\Pic(S_{\overline\bbQ}) \otimes_\bbZ \bbC}$
splits into irreducible components as
$\smash{\chi_\triv^{12} \oplus \chi_{\bbQ(\sqrt{-2})}}$.
Further~generators are formed by a line and two conics
in~$\Pb^2_\bbQ$,
the preimages of which split
in~$S$.
From~these altogether, one calculates~that
$$\Pic(S_{\overline\bbQ}) \otimes_\bbZ \bbC = \chi_\triv^{12} \oplus \chi_{\bbQ(\sqrt{-2})} \oplus \chi_{\bbQ(\sqrt{2})} \oplus \chi_{\bbQ(\sqrt{6})} \oplus \chi_{\bbQ(\sqrt{-6})}$$
and, consequently,
$\Delta_{\Pic}(S) = 1$.
The~splitting field of
$\smash{\Pic(S_{\overline\bbQ})}$
is~$\bbQ(i, \sqrt{2}, \sqrt{3})$.
On~the other hand, Algorithm~\ref{alg_delta} yields
$\Delta_{H^2}(S) = 3$,
so the jump character is given
by~$(\frac3\cdot)$.

If the effect of real multiplication added up with that of the transcendental character then this would imply
$\smash{\rk\Pic S_{\overline\bbF_{\!p}} > 18}$
for every
prime~$p$
such that
$\smash{(\frac3p) = -1}$,
a~contradiction.
\end{exo}

\section{Infinitely many rational curves}

It~has since long been conjectured that every
$K3$~surface
$S$
over an algebraically closed
field~$K$
contains infinitely many rational~curves. The problem has been settled 
only recently by X.~Chen, F.~Gounelas, and C.~Liedtke~\cite{CGL}.
Many particular cases had been known before, most notably, that of odd Picard rank (\cite{LL}, based on the ideas of~\cite{BHT}, cf.~\cite{Ben15}). Other sufficient conditions included those that
$S$
has infinitely many automorphisms, that
$S$
is elliptic~\cite{BT}, or that
$K$~is
of characteristic zero and
$S$
cannot be defined
over~$\overline\bbQ$~\cite[Theorem~3]{BHT}.

As an application of Theorem~\ref{K3_fl}, we show that the existence of infinitely many rational curves may be obtained rather easily in the situation that
$S$
is defined
over~$\overline\bbQ$,
the jump character is non-trivial and the surface is otherwise~generic. Our~result is as~follows.

\begin{theo}
\label{ratcrv}
Let\/~$K$
be a number field and\/
$S$
a\/
$K3$~surface
over\/~$K$.
Assume that\/
$\rk \Pic S_{\overline{K}}$
is even, that
$S_{\overline{K}}$
has neither real nor complex multiplication, and that\/
$\Delta_{H^2}(S)\Delta_{\Pic}(S)$
is a non-square
in\/~$K$.\smallskip

\noindent
Then\/
$S_{\overline{K}}$
contains infinitely many rational~curves.
\end{theo}

\begin{rem}
The transcendental part of
$\smash{T \subset H^2(S_\bbC, \bbQ)}$,
considered as a pure
\mbox{weight-$2$}
Hodge structure, has an endomorphism algebra
$\End_\Hs(T)$
that may only be a totally real field or a CM field~\cite[Theorem~1.6.(a) and Theorem~1.5.1]{Za}. Our~assumption concerning real and complex multiplication just means that
$\End_\Hs(T) = \bbQ$,
which is fulfilled as long as
$X$
is sufficiently~general.
\end{rem}

\begin{lem}
\label{nonsup}
Let\/~$K$
be a number field and\/
$S$
a\/~$K3$~surface
over\/~$K$.
Assume~that\/
$S_{\overline{K}}$
has neither real nor complex~multiplication. Then,~for every quadratic field
extension\/~$L/K$,
there are infinitely many inert primes\/
$\smash{\frakp \subset \calO_K}$
such that the
reduction\/~$\smash{S_{\overline\bbF_{\!\frakp}}}$
is~non-supersingular.\medskip

\noindent
{\bf Proof.}
{\em
The~case that
$\rk \Pic S_{\overline{K}} = 20$
is degenerate of the kind that
$T$
contains no
\mbox{$(1,1)$-classes}.
It~is known in this situation that
$S_{\overline{K}}$
automatically has complex multiplication \mbox{\cite[Theorem~4]{SI}}. We may therefore assume
that~$r := \dim T \geq 3$.

We~choose a
prime~$l$
and put
$\smash{T_l \subset H^2_\et(S_{\overline\bbQ}, \bbQ_l)}$
to be the transcendental part of
\mbox{$l$-adic}
cohomology. Then
$\smash{\GO(T_l, \langle.\,,.\rangle)_{\overline\bbQ_l}}$
is an algebraic group
over~$\bbQ_l$.
One~has
$\smash{\GO(T_l, \langle.\,,.\rangle)_{\overline\bbQ_l} \cong \GO_{r,\overline\bbQ_l}}$,
so that there are two connected components when
$r$~is
even, while the group is connected for
odd~$r$.

As~$\End_\Hs(T) = \bbQ$,
we know that the image of the canonical continuous representation
$\varrho_l\colon \Gal(\overline{K}/K) \to \GO(T_l, \langle.\,,.\rangle)$
is Zariski dense either
in~$\GO(T_l, \langle.\,,.\rangle)$
or in the neutral component
$\GO^0(T_l, \langle.\,,.\rangle)$.
Indeed,~this follows from the Mumford--Tate conjecture, proven by S.\,G.~Tankeev \cite{Ta90,Ta95}, together with Yu.\,G.~Zarhin's explicit description of the Mumford--Tate group in the case of a
$K3$~surface~\cite[Theorem~2.2.1]{Za}.

Now,~let us assume, to the contrary, that for all but finitely many inert
primes~$\frakp$,
the reduction
$\smash{S_{\overline\bbF_{\!\frakp}}}$
were~supersingular. We~put
\begin{align*}
M := \{&\frakp \subset \calO_K \text{ prime ideal} \mid \\[-1mm]
& \bbF_{\!\frakp} \text{ is a prime field,} \; \#\bbF_{\!\frakp} \neq l ,\;\frakp \text{ inert in } L, \;\frakp \text{ good for } S, \;S_{\overline\bbF_{\!\frakp}} \text{ supersingular}\} \, .
\end{align*}
Then~$M \subseteq I$,
for
$I$
the set of all inert primes, and the difference
$I \setminus M$
is of analytic density~zero. Indeed,~the prime ideals
$\frakp$
such that
$\bbF_{\!\frakp}$
is a prime field form a set of
density~$1$.

For~every prime ideal
$\frakp \subset \calO_K$,
we choose a geometric Frobenius automorphism
$\smash{\Frob_\frakp \in \Gal(\overline{K}/K)}$.
According~to the Chebotarev density theorem, the elements
$\smash{\sigma^{-1} \Frob_\frakp \sigma \in \Gal(\overline{K}/K)}$,
for
$\frakp \in M$
and
$\smash{\sigma \in \Gal(\overline{K}/K)}$,
are topologically dense in the non-trivial coset of
$\smash{\Gal(\overline{K}/K)}$
modulo
$\smash{\Gal(\overline{K}/L)}$.
Thus,~there are two elements
$\smash{\sigma_1, \sigma_2 \in \Gal(\overline{K}/K)}$
such that
$$\{\, \sigma_j \sigma^{-1} \Frob_\frakp \sigma \mid j = 1,2, \;\frakp \in M, \;\sigma \in \Gal(\overline{K}/K) \,\}$$
is dense in
$\smash{\Gal(\overline{K}/K)}$.

On the other hand, for
$\frakp \in M$
one has, due to supersingularity,
$\smash{p \;| \Tr \Frob_{\frakp, T_l}}$,
when writing
$p := \#\bbF_{\!\frakp}$.
Moreover,
$\smash{\det \Frob_{\frakp, T_l} = \pm p^r}$.
As~$|\Tr \Frob_{\frakp, T_l}| \leq rp$,
this shows~that
$$(\Tr \Frob_{\frakp, T_l})^r = \pm k^r \det \Frob_{\frakp, T_l} \,,$$
for
$-22 < -r \leq k \leq r < 22$.
Accordingly,~let
$C_k \subset \GO(T_l, \langle.\,,.\rangle)$
be the closed subscheme, defined by the equation
$(\Tr A)^r = \pm k^r \det A$,
and put
$C := \bigcup_{k=-r}^r C_k$.
Then~$\smash{C \subset \GO(T_l, \langle.\,,.\rangle)}$
is a closed subscheme and invariant under~conjugation.

As
$\smash{\GO(T_l, \langle.\,,.\rangle)_{\overline\bbQ_l} \cong \GO_{r,\overline\bbQ_l}}$,
for
$r \geq 3$,
it is easily seen that
$C$
cannot include a complete component of
$\smash{\GO(T_l, \langle.\,,.\rangle)_{\overline\bbQ_l}}$.
I.e.,~one has
$\dim C < \dim \GO(T_l, \langle.\,,.\rangle)$.
Thus,~the union
$\sigma_1 C \cup \sigma_2 C$
cannot be the whole group. Consequently,~the image of
$\smash{\Gal(\overline{K}/K) \to \GO(T_l, \langle.\,,.\rangle)}$
is Zariski dense neither in
$\GO(T_l, \langle.\,,.\rangle)$,
nor in
$\GO^0(T_l, \langle.\,,.\rangle)$,
a contradiction.
}
\eop
\end{lem}

\begin{prop}[Li--Liedtke]
\label{LL1}
Let\/~$K$
be a number field and\/
$S \subset \Pb^N_K$
a\/~$K3$~surface.
Assume that\/
$\Pic S = \Pic S_{\overline{K}}$
and that there is an infinite
set\/~$J$
of primes such that\/
$\rk \Pic S_{\overline\bbF_{\!\frakp}} > \rk \Pic S$
for\/~$\frakp \in J$.\smallskip

\noindent
Then~there exist a sequence without repetitions\/
$(\frakp_j)_{j\in\bbN}$
of primes
from\/~$J$
and a sequence\/
$(D_{\frakp_j})_{j\in\bbN}$
of rational curves\/
$\smash{D_{\frakp_j} \subset S_{\overline\bbF_{\!\frakp_j}}}$
such that the following two conditions are~satisfied.\smallskip

\noindent
The~class
$\smash{(D_{\frakp_j}) \in \Pic(S_{\overline\bbF_{\!\frakp_j}})}$
does not lie in the image
of\/~$\Pic S_{\overline{K}}$
under specialisation, for
any\/~$j$,
and\/~$\smash{\!\lim\limits_{j\to\infty}\! \deg D_{\frakp_j} = \infty}$.\medskip

\noindent
{\bf Proof.}
{\em
This is \cite[Proposition~4.2]{LL}.
\eop
}
\end{prop}

\begin{prop}[Li--Liedtke]
\label{LL_main}
Let\/~$K$
be a number field and\/
$S \subset \Pb^N_K$
a\/~$K3$~surface.
Assume that there exist a sequence without repetitions\/
$\smash{(\frakp_j)_{j\in\bbN}}$
of primes and a sequence\/
$(D_{\frakp_j})_{j\in\bbN}$
of rational curves\/
$\smash{D_{\frakp_j} \subset S_{\overline\bbF_{\!\frakp_j}}}$
\mbox{satisfying the following~conditions}.

\begin{iii}
\item
Each\/~$\smash{S_{\overline\bbF_{\!\frakp_j}}}$
is non-su\-persin\-gu\-lar,
\item
$(D_{\frakp_j})$
does not lie in the image
of the Picard
group\/~$\Pic S_{\overline{K}}$
of the generic fibre under~specialisation, for
any\/~$j$,
and\/
\item
$\smash{\lim\limits_{j\to\infty}\! \deg D_{\frakp_j} = \infty}$.
\end{iii}\smallskip

\noindent
Then, for every\/
$j \gg 0$,
there exists a rational curve\/
$D_j \subset S_{\overline{K}}$
such that its specialisation
to\/~$\smash{S_{\overline\bbF_{\!\frakp_j}}}$
is reducible, containing
$\smash{D_{\frakp_j}}$
as one of its~components. In~particular,
$\smash{\deg D_j > \deg D_{\frakp_j}}$.\medskip

\noindent
{\bf Proof.}
{\em
This is shown in the proof of~\cite[Theorem~4.3]{LL}.
}
\eop
\end{prop}

\begin{ttt}
{\bf Proof} of {\bf Theorem~\ref{ratcrv}.}
As~$\smash{\smash{\Delta_{H^2}(S)\Delta_{\Pic}(S)}}$
is a non-square
in~$K$,
the field
$\smash{L := K(\sqrt{\Delta_{H^2}(S)\Delta_{\Pic}(S)})}$
is indeed a quadratic~extension. By~Lemma~\ref{nonsup}, we have an infinite set
$J$
of inert primes such that
$\smash{S_{\overline\bbF_{\!\frakp}}}$
is~non-supersingular for every
$\frakp \in N$.
Moreover,
$\smash{\rk \Pic S_{\overline\bbF_{\!\frakp}} > \rk \Pic S_{\overline{K}}}$
according to Theorem~\ref{K3_fl}.b).

Let~now
$K' \supseteq K$
be the splitting field
of~$\Pic S_{\overline{K}}$.
For~each
$\frakp \in J$,
there is at least one prime
$\frakp' \subset \calO_{K'}$
lying
above~$\frakp$.
This~yields an infinite set
$J'$
of primes
in~$\calO_{K'}$,
to which Proposition~\ref{LL1} applies. It~provides a sequence
$(\frakp_j)_{j\in\bbN}$
of primes
in~$J'$
without repetitions and rational curves
$\smash{D_{\frakp_j} \subset S_{\overline\bbF_{\!\frakp_j}}}$,
not lying in the image
of~$\Pic S_{\overline{K}}$
under specialisation, such that
$\smash{\lim_{j\to\infty} \deg D_{\frakp_j} = \infty}$.
Knowing~this, Proposition~\ref{LL_main} yields a sequence
$\smash{(D_j)_{j\in\bbN}}$
of rational curves
on~$S_{\overline{K}}$
of degrees tending towards~infinity. This~completes the~proof.
\eop
\end{ttt}

\setlength\parindent{0mm}
\end{document}